\documentclass[a4paper]{article}

\font\cyr=wncyr10 % cyrillic font for Milnor's Lobachevski function.

\usepackage{graphicx}
\usepackage{amsmath}
\usepackage{amsthm}
\usepackage{amssymb}
\usepackage{subfigure}

\input{vacp-edge.sty}

\DeclareMathOperator{\im}{Im}
\DeclareMathOperator{\into}{\mathit{in}}
\DeclareMathOperator{\outof}{\mathit{ex}}

\newcommand{\e}[1]{e^{\displaystyle {#1}}}
\newcommand{\Li}[1][2]{\operatorname{Li}_{#1}}% dilogarithm
\newcommand{\Cl}[1][2]{\operatorname{Cl}_{#1}}% clausen's integral
\newcommand{\Sym}{\operatorname{Sym}}
\newcommand{\Alt}{\operatorname{Alt}}

\newcommand{\Seuc}{S_{\!\text{\it Euc}}}% the euclidean functional

\newcommand{\SeucH}{S_{\!\text{\it Euc}}^{\text{\it H}}}
\newcommand{\SeucHred}{\widehat{S}_{\!\text{\it Euc}}^{\text{\it H}}}
\newcommand{\Shyp}{S_{\!\text{\it hyp}}}% the hyperbolic functional

\newcommand{\ShypH}{S_{\!\text{\it hyp}}^{\text{\it H}}}
\newcommand{\ShypHred}{\widehat{S}_{\!\text{\it hyp}}^{\text{\it H}}}
\newenvironment{conditionlist}%
{\begin{list}{}%
{\setlength{\leftmargin}{\parindent}\setlength{\rightmargin}{\parindent}}}%
{\end{list}}

\theoremstyle{plain}

\newtheorem{theorem}{Theorem}
\newtheorem*{namedThm}{Theorem}
\newtheorem{proposition}{Proposition}
\newtheorem*{corollary}{Corollary}
\newtheorem*{flowThm}{Feasible Flow Theorem}

\theoremstyle{definition}

\newtheorem*{definition}{Definition}

\theoremstyle{remark}

\newtheorem*{remark}{Remark}

\begin{document}

\title{Variational principles for circle patterns and Koebe's theorem}
\author{Alexander I.~Bobenko\footnote{E--mail: {\tt bobenko@math.tu-berlin.de}}
\and 
Boris A.~Springborn\footnote{E--mail: {\tt springb@math.tu-berlin.de}}}

%\date{Version 1.21 of \today}
\date{Institut f\"ur Mathematik,
Technische Universit\"at Berlin, \\
Strasse des 17. Juni 136, 10623 Berlin, Germany}

\maketitle

%\tableofcontents\newpage

%---------------------------------------------------------------------

\section{Introduction}

\begin{figure}
\hfill
\subfigure[{}]{\label{fig:torus-combi}%
\includegraphics[width=0.45\textwidth]{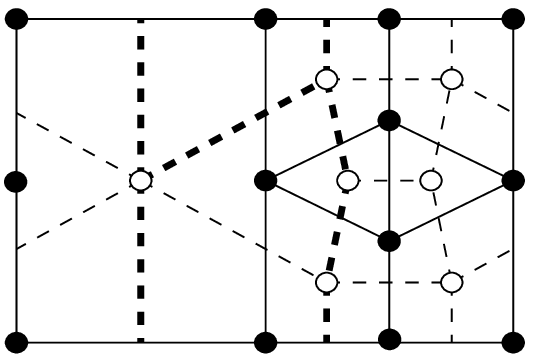}}
\hfill
\subfigure[{}]{\label{fig:torus-pattern}
\includegraphics[width=0.45\textwidth]{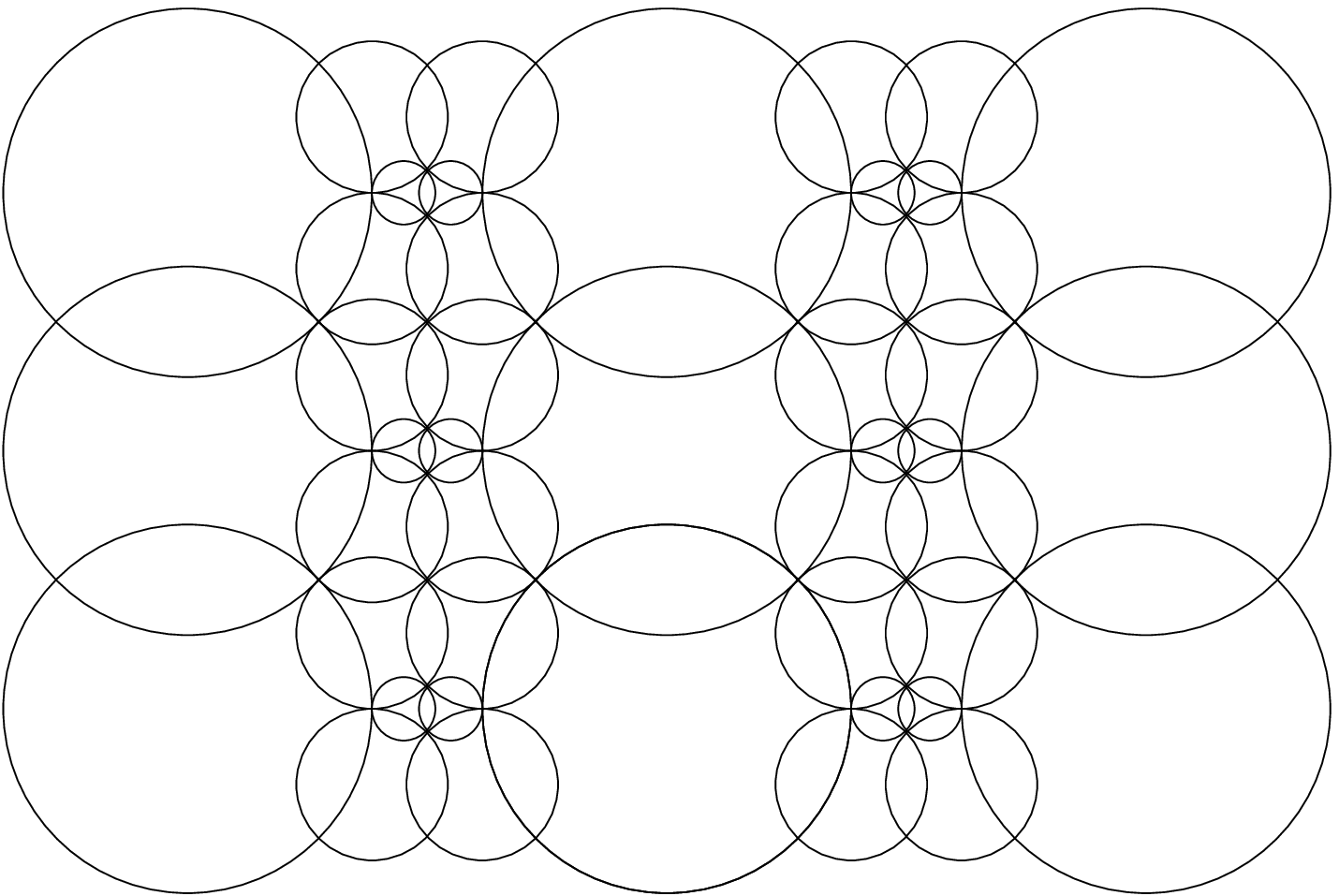}}
\hspace*{\fill}
\caption{A cellular decomposition of the torus and the corresponding circle
pattern with orthogonally intersecting circles.}
\label{fig:ex-torus}
\end{figure}

\begin{figure}
\hfill
\subfigure[{}]{\label{fig:disc-combi}%
\includegraphics[width=0.45\textwidth]{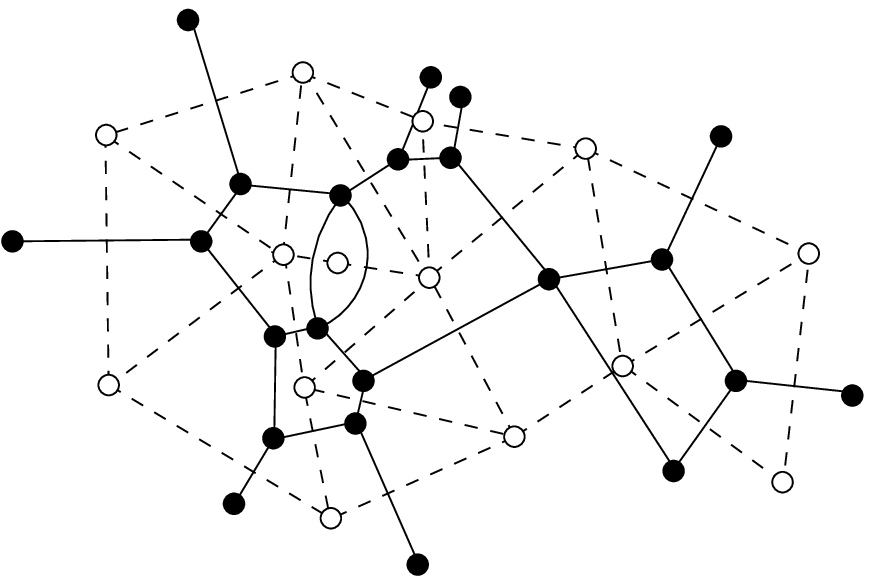}}
\hfill
\subfigure[{}]{\label{fig:disc-pattern}%
\includegraphics[width=0.45\textwidth]{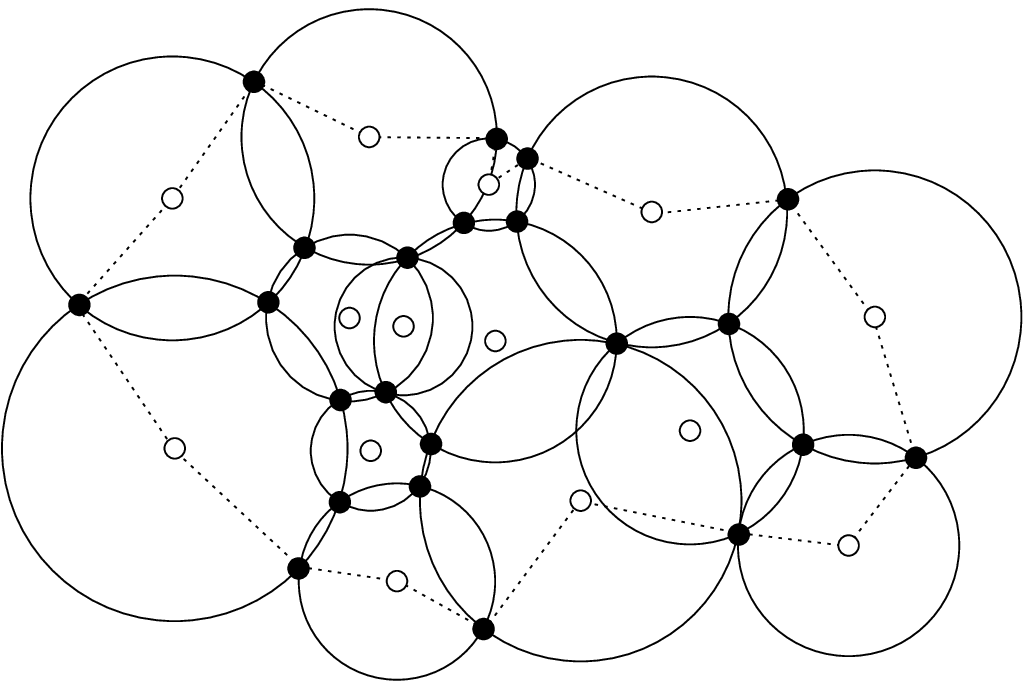}}
\hspace*{\fill}
\caption{A cellular decomposition of the disc and a corresponding circle
pattern with circles meeting at diverse angles.}
\label{fig:ex-disc}
\end{figure}

The subject of this paper is a special class of configurations, or patterns,
of intersecting circles in constant curvature surfaces. The combinatorial
aspect of such a pattern is described by a cellular decomposition of the
surface. The faces of the cellular decomposition correspond to circles and
the vertices correspond to points where circles intersect. (See
figures~\ref{fig:ex-torus} and~\ref{fig:ex-disc}.)

In the most general case that we consider, the surface may have cone-like
singularities in the centers of the circles and in the points of
intersection.

In particular, we treat the problem of constructing such circle patterns when
the cellular decomposition is given combinatorially, the intersection angles
of the circles are prescribed by a given function on the edges (and,
possibly, cone angles are prescribed by functions on the faces and vertices).
Using variational methods, and the `method of coherent angles,' we prove
existence and uniqueness results. The most fundamental one of these is
theorem \ref{thm:fundamental}. Similar results were obtained by Bowditch
\cite{Bowditch1991}, Garrett \cite{Garrett1992}, Rivin \cite{Rivin1999}, and
Leibon \cite{Leibon2001}. Our theorem treats circle patterns in surfaces
which may be flat or hyperbolic, closed surfaces or surfaces with boundary,
and there may be cone-like singularities in the centers of circles and in
their points of intersection. From it we deduce Rivin's theorem on ideal
hyperbolic polyhedra \cite{Rivin1996}, and theorem
\ref{thm:higherGenusRivin}, which is the analogous theorem for higher genus
surfaces. This higher genus analogue is new.

Our method is based on two new functionals---one for the Euclidean and one
for the hyperbolic case. We show how the functionals of Colin~de~Verdi{\`e}re
\cite{ColindeVerdiere1991}, Br{\"agger} \cite{Bragger1992}, and Rivin
\cite{Rivin1994} can be derived from ours. Leibon's functional
\cite{Leibon2001} seems to be related as well.

\section{A brief survey of relevant results}\label{sec:relevantResults}

A {\em circle packing}\/ is a configuration of disjoint discs which may touch
but not intersect. In 1936, Koebe published the following theorem about
circle packings in the sphere \cite{Koebe1936}.

%It was rediscovered by Thurston \cite{Thurston} as a
%consequence of Andreev's characterization of hyperbolic polyhedra
%\cite{Andreev1970a,Andreev1970b}.

\begin{namedThm}[Koebe]
For every triangulation of the sphere there is a packing of circles in the
sphere such that circles correspond to vertices and two circles touch if and
only if the corresponding vertices are adjacent. This circle pattern is
unique up to M{\"o}bius transformations of the sphere.
\end{namedThm}

Any graph without loops or double edges, embedded in the sphere, is
contained as subgraph in the $1$-skeleton of some triangulation. Hence, any
such graph is the touching graph of a circle pattern in the sphere. Of
course, the circle pattern is in general not unique up to M{\"o}bius
transformations.

There is, however, a way to generalize Koebe's theorem to cellular decompositions
of the sphere whose faces are not necessarily triangular {\em without giving
up the uniqueness part}. The right class of cellular decompositions to consider
is the class of cellular decompositions that are combinatorially equivalent to
convex polyhedra. By a theorem of Steinitz \cite{Steinitz1922,
SteinitzRademacher1934}, these {\em polytopal}\/ cellular decompositions are
characterized by the following two conditions. First, the cellular decomposition
is supposed to be {\em regular}, which means that that closure of each open
cell is a closed disc. For the $1$-cells (edges) this means that loops are
not allowed. For the $2$-cells this means that there are no identifications
on the boundary. The second condition is that two closed cells intersect in
{\em one}\/ closed cell, if at all. This forbids double edges. More
generally, it implies Steinitz' {\em Bedingung des Nich{\"u}bergreifens}\/:
If each of two faces is incident with each of two vertices, then there is an
edge which is incident with the two faces and the two vertices.

So given a polytopal cellular decomposition of the sphere, there is a
corresponding circle packing. Uniqueness is retained if an additional
condition is imposed on the packing: Consider the vertices incident with a
face of the cellular decomposition. The corresponding circles touch in points
which correspond to the edges of the face. It is required that these points
lie on a circle which intersects the original circles orthogonally. This is
the content of the following theorem.

\begin{theorem}\label{thm:orthoKoebe}
For every polytopal cellular decomposition of the sphere, there exists a pattern of
circles in the sphere with the following properties. There is a circle
corresponding to each face and to each vertex. The vertex circles form a
packing with two circles touching if and only if the corresponding vertices
are adjacent. Likewise, the face circles form a packing with circles touching
if and only if the corresponding faces are adjacent. For each edge, there is
a pair of touching vertex circles and a pair of touching face circles. These
pairs touch in the same point, intersecting each other orthogonally. 

This circle pattern is unique up to M{\"o}bius transformations.
\end{theorem}

If the cellular decomposition is a triangulation, the additional condition is
automatically fulfilled: Three mutually touching circles have three points of
contact. There is always a circle through them, and this circle intersects
the three original circles orthogonally. The above theorem is therefore
indeed a generalization of Koebe's theorem.

Furthermore, theorem \ref{thm:orthoKoebe} is equivalent to the following
theorem, which ``probably yields the best and most powerful version of
Steinitz' theorem \cite{Ziegler1995}.'' \footnote{The original Steinitz
  theorem only asserts the existence of a three dimensional convex polyhedron
  which is combinatorially equivalent to a polytopal cell decomposition of
  the sphere.}

\begin{theorem}\label{thm:superSteinitz}
For every polytopal cellular decomposition of the sphere, there is a combinatorially
equivalent polyhedron with edges tangent to a sphere. This polyhedron is
unique up to projective transformations which fix the sphere. 

There is a simultaneous realization of the dual polyhedron, such that
corresponding edges of the dual and the original polyhedron touch the sphere
in the same points and intersect orthogonally.
\end{theorem}

The history of theorems \ref{thm:orthoKoebe} and \ref{thm:superSteinitz} is
somewhat unclear. Sachs' survey \cite{Sachs1994} provides valuable
information. However, some of the articles in preparation which are cited
have never appeared in print. Thurston \cite{Thurston}, who revived the
interest in circle packings and circle patterns by pointing out the
connection to `discrete conformal maps' and the theory of $3$-manifolds,
devised a method to construct some kinds of circle patterns and thus prove
existence and uniqueness theorems. These methods were adapted by Brightwell
and Scheinerman \cite{BrightwellScheinerman1993}, who state and prove theorem
\ref{thm:orthoKoebe}. Schramm \cite{Schramm1992} proves a theorem that is
much more general than theorem \ref{thm:superSteinitz}. In it, the sphere is
replaced by an arbitrary convex body. The proof is not constructive.

\begin{figure}
\hfill
\subfigure[{}]{\label{fig:medial-a}%
\includegraphics[width=0.45\textwidth]{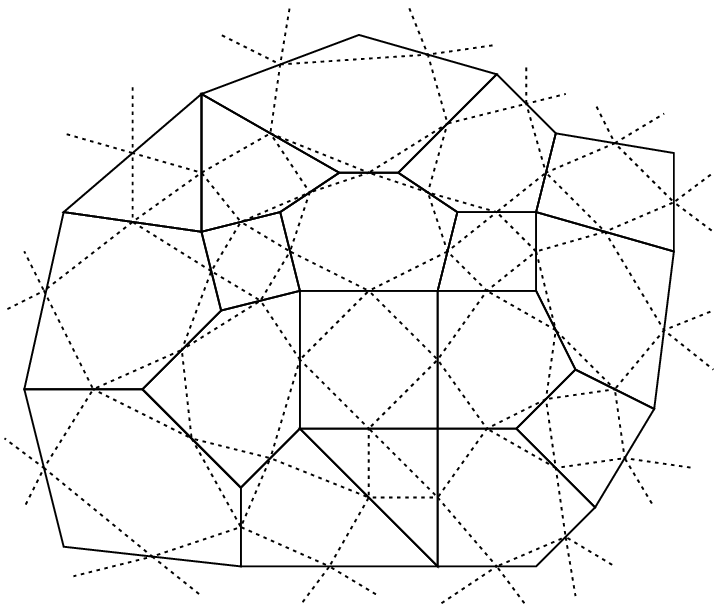}}
\hfill
\subfigure[{}]{\label{fig:medial-b}%
\includegraphics[width=0.45\textwidth]{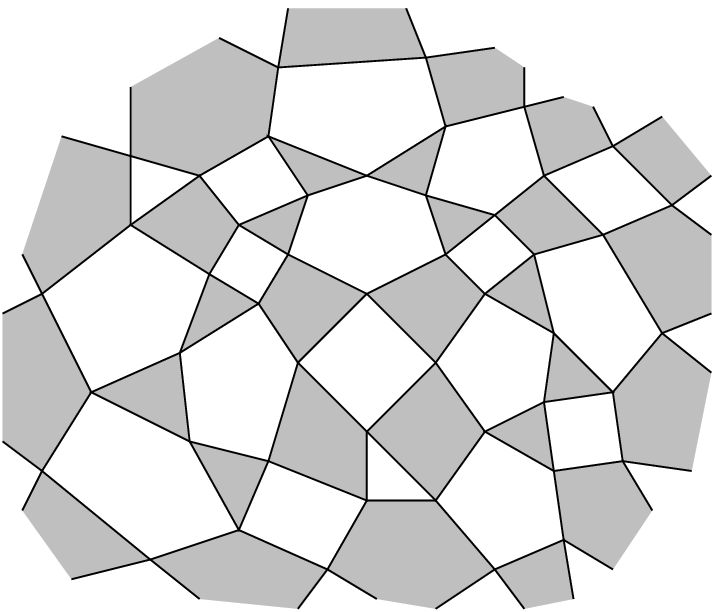}}
\hspace*{\fill}
\caption{A cellular decomposition (a) and its medial decomposition (b).}
\label{fig:medial}
\end{figure}

Theorem \ref{thm:orthoKoebe} may be generalized in two directions. First, one
may consider patterns of circles intersecting not orthogonally, but at
arbitrary angles. There is a correspondence between such circle patterns in
the sphere and polyhedra in hyperbolic space: In the Poincar{\'e} ball model,
hyperbolic planes are represented by spheres that intersect the unit sphere
orthogonally. Hence there is a correspondence between circles in the unit
sphere and hyperbolic planes. Furthermore, the intersection angle of two
circles equals the intersection angle of the corresponding
planes. Orientation preserving isometries of hyperbolic space correspond to
M{\"o}bius transformations of the sphere at infinity and vice versa. We
consider circle patterns that correspond to hyperbolic polyhedra with
vertices in the sphere at infinity, so called `ideal' hyperbolic
polyhedra. These are classified by the following theorem, which is due to
Rivin \cite{Rivin1996}. 
\begin{namedThm}[Rivin]
Let $\Sigma$ be a polytopal cellular decomposition of the sphere. Suppose an angle
$\theta_e$, with $0<\theta_e<\pi$, is assigned to each edge $e$. There exists
an ideal hyperbolic polyhedron which is combinatorially equivalent to
$\Sigma$, and with exterior dihedral angles $\theta_e$, if and only if the
following condition holds: For every cocycle $\gamma$ of edges,
$\sum_{e\in\gamma}\theta_e\geq 2\pi$, with equality if and only if $\gamma$
is the coboundary of a single vertex. 

This ideal hyperbolic polyhedron is unique up to isometry.
\end{namedThm}
Here, the {\em exterior intersection angle of two faces of a polyhedron} is
the complement to $\pi$ of the interior intersection angle, which is the
angle inside the polyhedron. The {\em exterior intersection angle of two
  circles}\/ is the angle $\theta_e$ of figure \ref{fig:kite}; it is the
complement to $\pi$ of the interior intersection angle $\theta_e^*$. A {\em
  cocycle}\/ is a cycle in the Poincar{\'e} dual cellular decomposition,
faces of which correspond to vertices of the original decomposition, vertices
to faces, and edges to edges. The concept dual to the boundary of a face is
the {\em coboundary}\/ of a vertex.
\begin{corollary}
Let $\Sigma$ be a polytopal cellular decomposition of the sphere, and suppose
every vertex has $n\geq 3$ edges. Suppose that every cocycle which is not the
coboundary of a single vertex is more than $n$ edges long. Then there exists,
uniquely up to M{\"o}bius transformations, a circle pattern on the sphere
which is combinatorially equivalent to $\Sigma$ and has exterior intersection
angles $2\pi/n$.
\end{corollary}
The case $n=4$ implies theorem \ref{thm:orthoKoebe}, and hence theorem
\ref{thm:superSteinitz}. Indeed, suppose $\Sigma$\/ is a polytopal cellular
decomposition of the sphere. In theorem \ref{thm:orthoKoebe}, circles
correspond to faces and vertices. To apply the corollary, consider the {\em
  medial}\/ cellular decomposition of $\Sigma$, whose faces correspond to the
faces and vertices of $\Sigma$ and whose vertices correspond to edges of
$\Sigma$. Figure \ref{fig:medial-a}
shows part of a cellular decomposition of the sphere and figure
\ref{fig:medial-b} its medial. The dotted lines in figure \ref{fig:medial-a}
represent the edges of the medial decomposition. In \ref{fig:medial-b}, the
faces of the medial decomposition which correspond to vertices in the
original decomposition are shaded. To deduce theorem \ref{thm:orthoKoebe}
from the corollary, one has to show that if we start with the polytopal
cellular decomposition $\Sigma$, its medial will be polytopal, have $4$-valent
vertices, and all cocycles which are not coboundaries of a single vertex will
be longer than $4$. In fact, this is easy to see. Suffice it to remark that
Steinitz' {\em Bedingung des Nicht{\"u}bergreifens}\/ for $\Sigma$\/
corresponds to the property of the medial that cocycles are longer than $4$.

The second direction in which theorem \ref{thm:orthoKoebe} may be generalized
is to consider not only circle patterns in the sphere but also in other
surfaces of constant curvature. We will consider circle patterns in closed
compact surfaces of arbitrary genus, as well as circle patterns in surfaces
with boundary. 

It is enough to consider only orientable surfaces. To deal with non-orientable
surfaces, consider first the orientable double cover; the uniqueness of the
circle patterns in question will imply that the pattern on the cover can be
projected down to the non-orientable surface.

By way of a further generalization, one may also consider circle patterns in
surfaces which may have cone-like singularities in the centers of
circles and in the intersection points. 

Colin~de~Verdi{\`e}re \cite{ColindeVerdiere1991} considers Euclidean and
hyperbolic circle packings, Br{\"a}gger \cite{Bragger1992} only the Euclidean
case.  Rivin \cite{Rivin1994} and Leibon \cite{Leibon2001} treat Euclidean
and hyperbolic circle patterns, respectively. They allow cone-like
singularities in the vertices.

In this paper, we present new variational principles for circle patterns in
Euclidean surfaces (with curvature $0$) and in hyperbolic surfaces (with
curvature $-1$). The functionals are functions of the radii of the
circles. Critical points correspond to circle patterns with prescribed
intersection angles and cone-like singularities. Circle patterns in the sphere
are treated by projecting them stereographically to the plane and then using
the Euclidean functional.

We give a complete proof of Rivin's theorem, and we prove an analogous
theorem for surfaces of higher genus, theorem
\ref{thm:higherGenusRivin}. Both are deduced from the more fundamental,
but less geometric theorem \ref{thm:fundamental}.

The idea to treat circle patterns and circle packings using a variational
principle is not new. Colin~de~Verdi{\`e}re \cite{ColindeVerdiere1991} was
the first to present a variational principle for circle packings, in which
the circles correspond to the vertices of a triangulation. He treats both the
Euclidean and the hyperbolic case. However, only the derivatives of the
functionals are given. This suffices to prove the existence and uniqueness of
circle packings.  Br{\"a}gger \cite{Bragger1992} remarks that, investigating
Colin~de~Verdi{\`e}re's functionals, one quickly arrives at the Lobachevski
function. But he does not present these functionals in closed form. Instead,
he derives a new variational principle.  Rivin \cite{Rivin1994}, and Leibon
\cite{Leibon2001} present variational principles for circle patterns. The
circles correspond to the faces of a triangulation. Rivin treats the
Euclidean, Leibon the hyperbolic case.

It was unclear how these functionals were related. We show how the
functionals of Colin~de~Verdi{\`e}re, Br{\"a}gger, and Rivin can be derived
from our functionals. We believe that Leibon's functional can be identified
with a functional we derive in section \ref{sec:other}.

Finally, we would like to add some remarks on the practical side of actually
constructing circles patterns. The functionals of Colin~de~Verdi{\`e}re are
not presented in closed form; only formulas for their derivatives are
given---and in a form which seems to have prevented their integration. The
functionals of Br{\"a}gger, Rivin, and Leibon are presented in closed form,
but their minimization has to be constrained to spaces of `coherent angle
systems', which makes their application more difficult. Our functionals are
given in closed form and there are no constraints. 

Thurston's proof of Koebe's theorem is also constructive.  ``[It is] based on
a practical algorithm for actually constructing patterns of circles. The idea
of the algorithm is to adjust, iteratively, the radii of the circles. A
change of any single radius affects most strongly the curvature at that
vertex, so this process converges reasonably well.''~\cite[p.~347]{Thurston}
Mohar shows that such an algorithm runs in polynomial time \cite{Mohar1993}.
Stephenson's program \textsf{circlepack} constructs circle packings by
Thurston's method. It was used by Dubejko and Stephenson for experiments in
discrete analytic function theory \cite{DubejkoStephenson1995a}. There are
also efforts to apply such circle packing algorithms in medical imaging
\cite{HurdalEtAl1999, Martindale2001}.

Thurston's algorithm can be interpreted as a special method to minimize our
functionals: Adjusting a single radius so that the neighboring circles fit
around it is equivalent to minimizing our functionals in one coordinate
direction.

%---------------------------------------------------------------------

\section{Main results}\label{sec:main}

After some preliminaries, we state the main results, theorem
\ref{thm:fundamental} and theorem \ref{thm:higherGenusRivin}. Theorem
\ref{thm:fundamental} is a fundamental (but somewhat technical) existence and
uniqueness theorem for circle patterns. Section
\ref{sec:mainproof} is devoted to the proof of theorem \ref{thm:fundamental}.
From it, we deduce (in section \ref{sec:sphereRivinHigherGenus}) Rivin's
theorem and its higher genus analogue, theorem \ref{thm:higherGenusRivin}.

Suppose a graph is embedded in a closed topological surface such that the
complement of the graph is the disjoint union of discs. This is called a {\em
  $2$-cell embedding} of the graph. Such a $2$-cell embedding of a graph
gives rise to a cellular decomposition of the surface. Two decompositions are
considered equivalent if one is mapped to the other by a homeomorphism of the
surface. We only consider the case of oriented surfaces. 
% In appendix \ref{app:cellular}, we elaborate on a combinatorial description
% of the resulting cellular decompositions. There, we also explain why we
% choose this rather broad class of decompositions.  For now, suffice it to say
% that
A {\em cellular surface} is a cellular decomposition arising from a $2$-cell
embedding, considered up to homeomorphism. Associated with a cellular surface
is a set $F$ of faces, a set $\vec{E}$ of oriented edges, and a set $V$ of
vertices. By $E$ we denote the set of unoriented edges. For each unoriented
edge in $E$, there is a pair of oppositely oriented edges in $\vec{E}$.

In the case of a surface with boundary, consider a graph embedded in the
surface, such that the graph meets the boundary only in vertices, and such
that the complement of the graph and the boundary is a disjoint union of open
discs. This gives rise to a cellular decomposition. As the edges of the
corresponding cellular surface we consider only the edges of the embedded
graph. This means that there are no boundary edges: Every edge has a face on
the left and on the right.  But there are faces whose boundary is not closed,
the {\em boundary faces}, and there are vertices whose coboundary is not
closed, {\em the boundary vertices.}

A convenient way to represent a cellular surface $\Sigma$ jointly with its
dual $\Sigma^*$ is the dual of the medial decomposition, $\Sigma_q$, which is
called a {\em quad-graph}, since all its faces are quadrilaterals. They
correspond to the edges of $\Sigma$. There are two kinds of vertices of
$\Sigma_q$, `white' vertices corresponding to faces of $\Sigma$ and `black'
vertices corresponding to vertices of $\Sigma$. Quad-graphs turn out to be
crucial for defining discrete integrable systems \cite{BobenkoSuris2002} (in
particular those related to circle patterns) and discrete Riemann surfaces
\cite{Mercat2001}.

\begin{figure}
  \centering
  \input{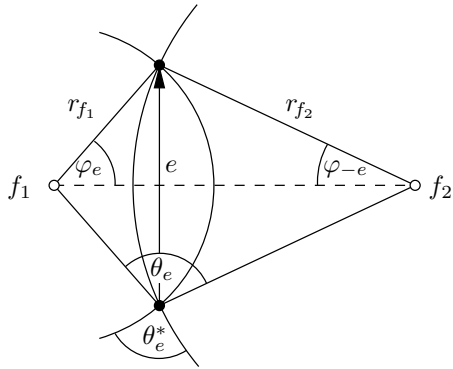}
  \caption{A Euclidean kite-shaped quadrilateral, $\theta_e+\theta^*_e=\pi$.}
  \label{fig:kite}
\end{figure}

Suppose we are given a positive radius $r_f$ for each face $f$ of $\Sigma$
and an angle $\theta_e\in(0,\pi)$ for each non-oriented edge of $\Sigma$. For
each face of $\Sigma_q$, construct a (Euclidean, spherical, or hyperbolic)
kites-shaped quadrilateral as in figure \ref{fig:kite}.
Because corresponding sides have the same length, these kites may be glued to
together to form a cellular decomposition equivalent to $\Sigma_q$. One
obtains not only a cellular decomposition of a topological manifold but a
{\em circle pattern} on a constant curvature surface with cone-like
singularities. The circles correspond to faces of $\Sigma$. The vertices of
$\Sigma$ correspond to points of intersection.

The cone-like singularities are in the centers of circles and in the points
corresponding to vertices.  The cone angle at an interior vertex $v\in
V_{\Sigma}$ is
\begin{equation}\label{eqn:vertexAngle}
\Theta_v=\sum\theta_e, 
\end{equation}
where the sum is taken over all oriented edges with initial vertex $v$. That
is, we interpret $\theta$ as a function on the oriented edges with $\theta_e
= \theta_{-e}$. (If an edge $e$ has $v$ as initial and terminal vertex, then
$\theta_e$ appears twice in the sum.)

The cone angle in the center of a circle corresponding to an interior face
$f$ is
\begin{equation}\label{eqn:bigPhi}
\Phi_f=2 \sum\varphi_e,
\end{equation}
where the sum is taken over all oriented edges in the oriented boundary of
$f$. An elementary geometric calculation yields for the angles 
\begin{equation}\label{eqn:phiOfREuc}
\varphi_e=
\frac{1}{2i}\log
\frac{r_{f_j}-r_{f_k}\, \e{-i\theta_e}}
     {r_{f_j}-r_{f_k}\, \e{i\theta_e}}
\end{equation}
in the Euclidean case, and 
\begin{multline}\label{eqn:phiOfRHyp}
\varphi_e=
\frac{1}{2i}
\left(
\log
\frac{\tanh\left(\frac{\displaystyle r_{f_j}}
                      {\displaystyle 2}
           \right)
     -\tanh\left(\frac{\displaystyle r_{f_k}}
                      {\displaystyle 2}
           \right) \e{-i\theta_e}}
     {\tanh\left(\frac{\displaystyle r_{f_j}}
                      {\displaystyle 2}
           \right)
     -\tanh\left(\frac{\displaystyle r_{f_k}}
                      {\displaystyle 2}
           \right) \e{i\theta_e}}  
\;-
\right.\\
\left.
\log
\frac{1-
     \tanh\left(\frac{\displaystyle r_{f_j}}
                      {\displaystyle 2}
           \right)
     \tanh\left(\frac{\displaystyle r_{f_k}}
                      {\displaystyle 2}
           \right) \e{-i\theta_e}}
     {1-
     \tanh\left(\frac{\displaystyle r_{f_j}}
                      {\displaystyle 2}
           \right)
     \tanh\left(\frac{\displaystyle r_{f_k}}
                      {\displaystyle 2}
           \right) \e{i\theta_e}}
\right)
\end{multline}
in the hyperbolic case. Here, $f_j$ and $f_k$ are the faces to the
left and right of $e$.

For boundary faces and vertices, $\Phi_f$ and $\Theta_v$, given by the same
formulas, are the interior angles of the boundary polygon of the
quadrilateral decomposition (the dotted polygon in figure
\ref{fig:disc-pattern}).

The idea is to prescribe only the intersection angles $\theta_e$ and the cone
(or boundary) angles $\Phi_f$ in the faces, and to solve for the radii. The
cone (or boundary) angles $\Theta_v$ in the vertices are determined by the
intersection angles and equation (\ref{eqn:vertexAngle}). We will prove the
following theorem.
\begin{theorem}\label{thm:fundamental}
Let $\Sigma$ be an oriented cellular surface, let $\theta^*\in(0,\pi)^E$ be a
function on the non-oriented edges, and let $\Phi\in(0,\infty)^{F_{\Sigma}}$
be a function on the faces.

A Euclidean circle pattern which is combinatorially equivalent to
$\Sigma$, which has interior intersection angles $\theta^*$, and which has
cone (or boundary) angles $\Phi_f$ in the centers of the circles exists, if
and only if the following two conditions are satisfied.
\begin{conditionlist}
\item{(i)}
\begin{equation}\label{eqn:eqCond}
\sum_{f\in F}\Phi(f)=\sum_{e\in E}2\theta^*(e)
\end{equation}
\item{(ii)}
If $F'$ is a nonempty subset of the face set $F$, $F'\not=F$, and 
$E'$ is the set of all edges incident with any face in $F'$, then
\begin{equation}\label{eqn:ineqCond}
\sum_{f\in F'}\Phi(f)<\sum_{e\in E'}2\theta^*(e).
\end{equation}
\end{conditionlist}
The Euclidean circle pattern is unique up to similarity.

A corresponding hyperbolic circle pattern in a surface of constant curvature
$-1$ with cone-like singularities exists, if and only if inequality
(\ref{eqn:ineqCond}) holds for all nonempty subsets $F'\subset F$, including
$F$ itself. The hyperbolic circle pattern is unique up to isometry.
\end{theorem}

This theorem is proved in section \ref{sec:mainproof}.

Equation (\ref{eqn:eqCond}), which is necessary for the existence of
a Euclidean circle pattern, is equivalent to the Gauss-Bonnet formula
\begin{equation*}
\sum_{f\in F} K_f + \sum_{v\in V} K_v = 2\pi\chi,
\end{equation*}
where $\chi$ is the Euler characteristic of $\Sigma$,
\begin{equation*}
K_f=\left\{
\begin{array}{ll}
2\pi-\Phi_f & \text{for interior faces} \\
 \pi-\Phi_f & \text{for boundary faces}
\end{array}
\right.,
\end{equation*}
and
\begin{equation*}
K_v=\left\{
\begin{array}{ll}
2\pi-\Theta_v & \text{for interior vertices} \\
\pi-\Theta_v & \text{for boundary vertices}
\end{array}
\right. .
\end{equation*}
Note that the formula $\chi=|F|-|E|+|V|$ does not hold for cellular surfaces
with boundary faces and boundary vertices. To get a topological invariant,
one must either add boundary edges first, or count the boundary faces and
vertices only half.

So far, we did not consider circle patterns in the sphere. The metric
approach, focusing on the radii of circles, is less useful in this case. The
reason is this: The only circle preserving transformations of the hyperbolic
plane are isometries, those of the Euclidean plane are similarities. The
M{\"o}bius transformations, which are the circle preserving transformations
of the sphere, distort the metric in a more complicated way.

However, spherical circle patterns {\em without cone-like singularities}\/ may
be treated by projecting them stereographically to the plane. This is
described in detail in section \ref{sec:sphereRivinHigherGenus}, where, in
particular, we give a proof of Rivin's theorem.

Searching for a generalization of Rivin's theorem to surfaces of higher
genus, one may conjecture the condition to be that the angle sum for a
coboundary is at least $2\pi$, and equal to $2\pi$ only for coboundaries of a
single vertex. This is not the right condition, however. Firstly, one needs
to consider only the coboundaries of simply connected domains. Secondly, it
is not enough to consider coboundaries. For example, if we cut the torus of
figure \ref{fig:torus-combi} along the dual edges represented by the heavy
dashed lines, it is separated into two pieces; one of these pieces is a
disc. While it is indeed necessary that the angle sum over the edges we cut
along is strictly greater than $2\pi$, these edges do not form a boundary. To
formulate the right condition, we consider cutting the dual cellular surface
along some of its edges. Of the pieces into which this separates the surface,
only the discs have to be considered. The angle sum over their boundary
has to be at least $2\pi$. What counts is the boundary of the pieces after
cutting, not their boundary as regions in the surface.

\begin{theorem}\label{thm:higherGenusRivin}
Let $\Sigma$ be a cellular surface of genus $g>0$, and let $\Sigma^*$ be the
dual cellular surface. Suppose $\theta\in(0,\pi)^E$ is a function on the
edges of $\Sigma$, which sums to $2\pi$ around each vertex. We identify the
edges of $\Sigma$ with the edges of $\Sigma^*$ and consider $\theta$ also as
function on the edges of $\Sigma^*$.

There exists a circle pattern in a surface of constant curvature $-1$ which
is combinatorially equivalent to $\Sigma$ and has exterior intersection
angles $\theta$, if and only if the following condition is satisfied.
\begin{conditionlist}
\item{} Suppose we cut $\Sigma^*$ along a subset of its edges,
obtaining one or more pieces. For any piece that is a disc, the sum of
$\theta$ over the edges in its boundary has to be at least $2\pi$, with
equality if and only if the piece consists of only one face of $\Sigma^*$
(vertex of $\Sigma$).
\end{conditionlist}
Both the constant curvature surface and the circle pattern are determined
up to isometry if $g>1$ and up to similarity if $g=1$. 
\end{theorem}

The proof is contained in section \ref{sec:sphereRivinHigherGenus}.

%---------------------------------------------------------------------

\section{The proof of theorem \ref{thm:fundamental}}\label{sec:mainproof}

In this section, we prove theorem \ref{thm:fundamental}. Our proof is based
on a variational principle for circle patterns. The corresponding
functionals, one for the Euclidean and one for the hyperbolic case, are
presented in section \ref{sec:funConv}. They are functions of the radii. Of
course, they also depend on the cellular surface; and the prescribed
intersection angles and cone angles enter as parameters. A radius function
which is a critical point yields a circle pattern with the prescribed
intersection angles and cone angles.  The functionals are strictly convex,
and this implies the uniqueness claims of theorem \ref{thm:fundamental}.

Proving the existence claims is a little more involved. One has to estimate
the functionals from below to show that they become large if some radii tend
to zero or infinity. In section \ref{sec:cas}, we define `coherent angle
systems', which are functions on the directed edges satisfying some linear
equations and inequalities. We show that the functionals have minima if and
only if a coherent angle system exists. This reduces the problem to the
feasibility problem of linear programming, or equivalently, to the feasible
flow problem of network theory.

In section \ref{sec:flowthm}, applying the feasible flow theorem of network
theory, we show that the conditions of theorem \ref{thm:fundamental} are
necessary and sufficient for the existence of a coherent angle system. 

%---------------------------------------------------------------------

\subsection{The Functionals. Convexity and the uniqueness of circle patterns}
\label{sec:funConv}

Instead of the radii $r$, it is more convenient to use transformed
variables $\rho$, the transformation being
\begin{equation*}
\rho = \left\{ 
\begin{array}{ll}
\log r           & \text{in the Euclidean case} \\
\log\tanh(r/2)   & \text{in the hyperbolic case}
\end{array}
\right..
\end{equation*}
In the Euclidean case, positive radii correspond to real $\rho$. In the
hyperbolic case case, positive radii correspond to negative $\rho$. 

The problem is to find a function $\rho: F\rightarrow \mathbb{R}$ (which must
be negative in the hyperbolic case) satisfying equation (\ref{eqn:bigPhi})
with
\begin{equation}\label{eqn:phiOfRhoEuc}
\varphi_e=
\frac{1}{2i}\log\frac{1-\e{\rho_{f_k}-\rho_{f_j}-i\theta_e}}
                     {1-\e{\rho_{f_k}-\rho_{f_j}+i\theta_e}},
\end{equation}
in the Euclidean case and
\begin{equation}\label{eqn:phiOfRhoHyp}
\varphi_e=
\frac{1}{2i}\log\frac{1-\e{\rho_{f_k}-\rho_{f_j}-i\theta_e}}
                     {1-\e{\rho_{f_k}-\rho_{f_j}+i\theta_e}}
-
\frac{1}{2i}\log\frac{1-\e{\rho_{f_k}+\rho_{f_j}-i\theta_e}}
                     {1-\e{\rho_{f_k}+\rho_{f_j}+i\theta_e}}.
\end{equation}
in the hyperbolic case. Here, $f_j$ and $f_k$ are the faces to the left and
right of $e$, respectively. (Equations (\ref{eqn:phiOfRhoEuc}) and
(\ref{eqn:phiOfRhoHyp}) are equivalent to equations (\ref{eqn:phiOfREuc}) and
(\ref{eqn:phiOfRHyp}).)

We turn this into a variational problem.  The corresponding functionals are
defined in terms of the imaginary part of the dilogarithm function $\Li$. In
the appendix, all requisite facts about this function and the
related Clausen's integral $\Cl$ are summarized.

\begin{definition}\label{def:functionals}
The {\em Euclidean circle pattern functional} is
\begin{equation}\label{eqn:SEuc}
\begin{split}
\Seuc(\rho)=
\sum_{f_j\circ\edge\circ f_k} \Big( &
 \im\Li\big(\e{\rho_{f_k}-\rho_{f_j}+i\theta_e}\big)
+\im\Li\big(\e{\rho_{f_j}-\rho_{f_k}+i\theta_e}\big)\\
&-\theta^*_e\big(\rho_{f_j}+\rho_{f_k}\big)\Big)
+\sum_{\circ f}\Phi_f \rho_f.
\end{split}
\end{equation}
The first sum is taken over all edges $e$, and $f_j$ and $f_k$ are the faces
on either side of $e$. The summands are symmetric in $f_j$ and $f_k$, so it
does not matter which face is considered as $f_j$ and which as $f_k$. The
second sum is taken over all faces $f$.

The {\em hyperbolic circle pattern functional} is, using the same notation, 
\begin{equation}\label{eqn:Shyp}
\begin{split}
\Shyp(\rho)= \sum_{f_j\circ\edge\circ f_k} \Big( &
  \im\Li\big(\e{\rho_{f_k}-\rho_{f_j}+i\theta_e}\big)
 +\im\Li\big(\e{\rho_{f_j}-\rho_{f_k}+i\theta_e}\big)\\
&+\im\Li\big(\e{\rho_{f_j}+\rho_{f_k}+i\theta_e}\big)
 +\im\Li\big(\e{-\rho_{f_j}-\rho_{f_k}+i\theta_e}\big)
   \Big)\\
+\sum_{\circ f} &\Phi_f\rho_f.
\end{split}
\end{equation}
\end{definition}

\begin{proposition}
The critical points of the functionals $\Seuc$ and $\Shyp$ correspond circle
patterns with intersection angles $\theta$ and cone angles $\Phi$ in the
centers of circles. In particular, if $\rho$ is a critical point of $\Shyp$
then all $\rho_f$ are negative.
\end{proposition}
\begin{proof}
The functionals are of the form
\begin{equation}\label{eqn:SForm}
S(\rho)=\sum_{f_j\circ\edge\circ f_k} 
\big(F_-(\rho_{f_j}-\rho_{f_k}) + F_+(\rho_{f_j}+\rho_{f_k}) \big)
+ \sum_{\circ f} G(\rho_f),
\end{equation}
where $F_-$ is an even function. The partial derivative with respect to a
face $\rho_{f_j}$ is  
\begin{equation*}
\frac{\partial S}{\partial\rho_{f_j}} =
\sum_{f_j\circ\edgeup\circ f_k} \big(F_-'(\rho_{f_j}-\rho_{f_k}) 
                                   + F_+'(\rho_{f_j}+\rho_{f_k})\big)
+ G'(\rho_{f_j}),
\end{equation*}
where the sum is taken over all oriented edges in the boundary of $f_j$.
Since, from equation (\ref{eqn:ImLiIntegral}) of the appendix,
\begin{equation*}
\frac{d}{dx}\,\im\Li(\e{x+i\theta})=\frac{1}{2i}\log\frac{1-\e{x-i\theta}}
                                                         {1-\e{x+i\theta}},
\end{equation*}
one finds that
\begin{equation}\label{eqn:partialS}
\frac{\partial S}{\partial \rho_f}=-2\big({\textstyle\sum}\varphi_e\big) 
+ \Phi_f,
\end{equation}
where the sum is taken over all directed edges in the boundary of $f$. 

It is left to show that if $\rho$ is a critical point of $\Shyp$, then all
$\rho_f$ are negative. The function
\begin{equation*}
f(x)=\frac{1}{2i}\log\frac{1-\e{x-i\theta_e}}
                          {1-\e{x+i\theta_e}}
\end{equation*}
is strictly increasing for $0<\theta<\pi$. Hence, if, in equation
(\ref{eqn:phiOfRhoHyp}), $\rho_{f_j}$ is nonnegative, then $\varphi_e$ is
non-positive. Therefore, if some $\rho_f$ is nonnegative, then, from equation
(\ref{eqn:partialS}),
\begin{equation}\label{eqn:dShypIneq}
\frac{\partial \Shyp}{\partial \rho_f}\geq\Phi(f)>0, 
\quad\text{if}\quad
\rho_f\geq 0.
\end{equation}
\end{proof}
 
Rescaling all radii by a factor $h$ corresponds to the transformation
$\rho\mapsto \rho+h\,1_F$, where $1_F$ is the function that is constantly $1$
on $F$. Note that
\begin{equation*}
\Seuc(\rho+h\,1_F)=\Seuc(\rho)
+h \Big( \sum_{f\in F} \Phi_f - 2 \sum_{e\in E} \theta^*_e \Big).
\end{equation*}
The functional can have a critical point only if the coefficient of $h$
vanishes. In this case, the functional is scale-invariant, and one may
restrict minimization to the subspace  
\begin{equation}\label{eqn:U}
U=\{\rho\in\mathbb{R}^F|\sum_{f\in F}\rho_f=0\}.
\end{equation}
\begin{proposition}
The hyperbolic functional is strictly convex. The Euclidean functional is
strictly convex on the subspace $U$. This implies, in the hyperbolic case,
the uniqueness of circle patterns with prescribed intersection angles and cone
angles, and, in the Euclidean case, uniqueness up to scale.
\end{proposition}
\begin{proof}
Taking second derivatives on both sides of equation (\ref{eqn:SForm}), one
obtains
\begin{equation*}
\begin{split}
S''=
\sum_{f_j\circ\edge\circ f_k} &\Big(
  F_-''(\rho_{f_j}-\rho_{f_k})\,(d\rho_{f_j}-d\rho_{f_k})^2 
+ F_+''(\rho_{f_j}+\rho_{f_k})\,(d\rho_{f_j}+d\rho_{f_k})^2 
\Big)\\
+ \sum_{\circ f} & G''(\rho_f)\,d\rho_f^2.
\end{split}
\end{equation*}
Thus, by a straightforward calculation, 
\begin{equation*}
\Seuc''=\sum_{f_j\circ\edge\circ f_k}
\frac{\sin\theta_e}
     {\cosh(\rho_{f_k}-\rho_{f_j})-\cos\theta_e}
\,(d\rho_{f_k}-d\rho_{f_j})^2
\end{equation*}
and
\begin{multline*}
\Shyp''=\sum_{f_j\circ\edge\circ f_k}\left(
\frac{\sin\theta_e}
     {\cosh(\rho_{f_k}-\rho_{f_j})-\cos\theta_e}
\,(d\rho_{f_j}-d\rho_{f_k}\big)^2 +
\right.
\\
\left.
\frac{\sin\theta_e}
     {\cosh(\rho_{f_j}+\rho_{f_k})-\cos\theta_e}
\,(d\rho_{f_j}+d\rho_{f_k}\big)^2
\right).
\end{multline*}
These second derivatives are positive definite on $U$ and $\mathbb{R}^F$,
respectively.
\end{proof}

\begin{remark}
Formulas (\ref{eqn:SEuc}) and (\ref{eqn:Shyp}) for $\Seuc$ and $\Shyp$
involve the imaginary part of the dilogarithm of a complex argument. The
formulas below involve only real functions. Since, moreover, Clausen's
integral is periodic, these formulas are well suited for the practical
evaluation of the functionals.

For $\rho\in\mathbb{R}^F$, let $p$ and $s$ be the functions on oriented edges
defined by
\begin{equation}\label{eqn:pOfRho}
p_e=2\arctan\left(
\tan\left(\frac{\theta^*_e}{2}\right)
\tanh\left(\frac{\rho_{f_k}-\rho_{f_j}}{2}\right)
\right)
\end{equation}
and
\begin{equation}\label{eqn:sOfRho}
s_e=2\arctan\left(
\tan\left(\frac{\theta^*_e}{2}\right)
\tanh\left(\frac{\rho_{f_k}+\rho_{f_j}}{2}\right)
\right),
\end{equation}
where $f_j$ and $f_k$ are the faces to the left and right of $e$. Note that
$s_e=s_{-e}$ (so that $s$ may be considered as a function on non-oriented
edges), but $p_e=-p_{-e}$.

The following formulas hold:
\begin{multline*}
\Seuc(\rho)=
\sum_{f_j\circ\edge\circ f_k}\Big(
  p_e(\rho_{f_k}-\rho_{f_j})
  +\Cl(\theta^*_e+p_e)
  +\Cl(\theta^*_e-p_e)
  -\Cl(2\theta^*_e)\\
  -\theta^*_e(\rho_{f_j}+\rho_{f_k})\Big)
+ \sum_{\circ f} \Phi_f\rho_f,
\end{multline*}
and
\begin{multline*}
  \Shyp(\rho)=
  \sum_{f_j\circ\edge\circ f_k} \Big( 
  p_e\big(\rho_{f_k}-\rho_{f_j}\big)
  +\Cl\big(\theta^*_e+p_e\big)+\Cl\big(\theta^*_e-p_e\big)\\
  +s_e\big(\rho_{f_j}+\rho_{f_k}\big)
  +\Cl\big(\theta^*_e+s_e\big)+\Cl\big(\theta^*_e-s_e\big)\\
  -2\Cl\big(2\theta^*_e\big) \Big)
  +\sum_{\circ f} \Phi_f\rho_f.
\end{multline*}
The first sum in each equation is taken over non-oriented edges, or rather,
over one oriented representative $e$ for each non-oriented edge. Since both
$p_e$ and $(f_k-f_j)$ change sign when $-e$ is chosen instead of $e$, the
summands are independent of which choice is made.

These formulas can be derived using equation (\ref{eqn:LiClp}) of
the appendix and some simple properties of Clausen's integral,
which are also found there.

The hyperbolic functional will reappear in this form in equation
(\ref{eqn:ShypH}) of section \ref{sec:other}.
\end{remark}

%---------------------------------------------------------------------

\subsection{Coherent angle systems and the existence of circle patterns}
\label{sec:cas}

In section \ref{sec:funConv}, the uniqueness of a circle pattern was deduced
from the convexity of the functionals. The rest of section
\ref{sec:mainproof} is devoted to the existence part of theorem
\ref{thm:fundamental}. To establish that the functionals have a minimum, we
show that they tend to infinity as $\rho$ tends to infinity in
$\mathbb{R}$ or $U$, respectively. 

The main difficulty in estimating the functionals from below is to somehow
merge the sum over edges with the sum over faces in equations
(\ref{eqn:SEuc}) and (\ref{eqn:Shyp}). This is achieved with the help of a so
called `coherent angle system'. Therefore, in this section we prove that
the functionals have minima if and only if coherent angle systems exist. In
section \ref{sec:flowthm}, we show that conditions of theorem
\ref{thm:fundamental} are necessary and sufficient for the existence of a
coherent angle systems.

\begin{definition}
A {\em coherent angle system} is a function $\varphi\in\mathbb{R}^{\vec{E}}$
on the oriented edges satisfying, in the Euclidean case, 
\begin{enumerate}
\item\label{cas1}
for all oriented edges $e\in\vec{E}$, 
$$\varphi_e>0\quad\text{and}\quad\varphi_e+\varphi_{-e}=\theta^*_e$$ 
\item\label{cas2}
for all oriented faces $f\in F$,
$$ \sum 2\varphi_e=\Phi_f, $$
where the sum is taken over all oriented edges in the oriented boundary of
$f$.
\end{enumerate}
In the hyperbolic case, condition \ref{cas1} has to be changed to 
\begin{enumerate}
\item[\ref{cas1}$'$.] 
for all oriented edges $e\in\vec{E}$, 
$$\varphi_e>0\quad\text{and}\quad\varphi_e+\varphi_{-e}<\theta^*_e.$$ 
\end{enumerate}
\end{definition}

\begin{proposition}
The functional $\Seuc$ ($\Shyp$) has a critical point, if and only if a
Euclidean (hyperbolic) coherent angle system exists.
\end{proposition}
\begin{proof}
If the functional $\Seuc$ ($\Shyp$) has a critical point $\rho$, then
equation (\ref{eqn:phiOfRhoEuc}) (equation (\ref{eqn:phiOfRhoHyp})) defines a
coherent angle system. It is left to show that, conversely, the existence of
a coherent angle system implies the existence of a critical point.

Consider first the Euclidean case. Suppose a Euclidean coherent angle system
$\varphi$ exists. One easily deduces equation (\ref{eqn:eqCond}). Hence, the
functional $\Seuc$ is scale invariant. We will show that
$\Seuc(\rho)\rightarrow\infty$ if $\rho\rightarrow\infty$ in the subspace $U$
defined in equation (\ref{eqn:U}). More precisely, we will show that for
$\rho \in U$,
\begin{equation}\label{eqn:SEucEstimate}
\Seuc(\rho)>2\min_{e\in\vec{E}}\varphi_e\,\max_{f\in F}\,|\rho_f|.
\end{equation}
The functional $\Seuc$ must therefore attain a minimum, which is a critical
point.

For $x\in\mathbb{R}$ and $0<\theta<\pi$,
\begin{equation*}
\im\Li(\e{x+i\theta})+\im\Li(\e{-x+i\theta})>(\pi-\theta)\,|x|,
\end{equation*}
and hence 
\begin{equation*}
\Seuc(\rho)>-2\sum_{e\in E}
\theta^*_e\min\big(\rho_{f_k},\rho_{f_j}\big) 
+ \sum_{f\in F}\Phi_f\rho_f,
\end{equation*}
Throughout this section, `$\sum_{e\in E}$' is meant to indicate that the sum
is taken over one representative from each pair of oppositely oriented edges;
it does not matter which. The faces on either side of $e$ are $f_j$ and
$f_k$.

Now we use the coherent angle system to merge the two sums. Because
\begin{equation*}
\sum_{f\in F}\Phi_f\rho_f 
= 2\sum_{e\in E} (\varphi_e\,\rho_{f_j} + \varphi_{-e}\,\rho_{f_k}), 
\end{equation*}
one obtains
\begin{equation*}
\Seuc(\rho)>2\sum_{e\in E}
\min\big(\varphi_e,\varphi_{-e}\big)\,
\big|\rho_{f_k}-\rho_{f_j}\big|.
\end{equation*}
Since we assume the cellular surface to be connected, we get
\begin{equation*}
\Seuc(\rho)>2\min_{e\in\vec{E}}\varphi_e\,
\big(\max_{f\in F}\rho_f-\min_{f\in F}\rho_f\big),
\end{equation*}
and from this the estimate (\ref{eqn:SEucEstimate}).

The hyperbolic case is similar. One shows that, if all $\rho_f<0$,
\begin{equation*}%\label{eqn:ShypEstimate}
\Shyp(\rho)>2
\min_{e\in E}\big|\varphi_e+\varphi_{-e}-\theta^*_e\big|\;
\max_{f\in F}\big|\rho_f\big|.
\end{equation*}
Since it was shown in section \ref{sec:funConv} that
$\partial\Shyp/\partial\rho_f>0$ if $\rho_f\geq0$, this implies the existence
of a minimum.
\end{proof}

%---------------------------------------------------------------------

\subsection{Conditions for the existence of coherent angle systems}
\label{sec:flowthm}

With this section we complete the proof of theorem \ref{thm:fundamental}. All
that is left to show is the following proposition. 
\begin{proposition}
A coherent angle system exists if and only if the conditions of theorem
\ref{thm:fundamental} hold.
\end{proposition}
\begin{proof}
It is easy to see that these conditions are necessary. To prove that they are
sufficient, we apply the feasible flow theorem of network theory. Let $(N,X)$
be a network (i.e.\ a directed graph), where $N$ is the set of nodes and $X$
is the set of branches. For any subset $N'\subset N$ let $\outof(N')$ be
the set of branches having their initial node in $N'$ but not their terminal
node. Let $\into(N')$ be the set of branches having their terminal node in
$N'$ but not their initial node. Assume that there is a lower capacity bound
$a_x$ and an upper capacity bound $b_x$ associated with each branch $x$,
with $-\infty\leq a_x\leq b_x\leq \infty$. 

\begin{definition}
A {\em feasible flow} is a function
$\varphi\in\mathbb{R}^X$, such that Kirchoff's current law is
satisfied, i.e., for each $n\in N$,
\begin{equation*}
\sum_{x\in \outof(\{n\})} \varphi_x = \sum_{x\in \into(\{n\})} \varphi_x,
\end{equation*}
and $a_x\leq\varphi_x\leq b_x$ for all branches $x$.
\end{definition}

\begin{flowThm}
A feasible flow exists, if and only if for every nonempty subset $N'\subset
N$ of nodes with $N'\not=N$, 
\begin{equation*}
\sum_{x\in \outof(N')} b_x \geq \sum_{x\in \into(N')} a_x.
\end{equation*}
\end{flowThm}

A proof is given by Ford and Fulkerson \cite[Ch.~II,
\S3]{FordFulkerson1962}. (Ford and Fulkerson assume the capacity bounds to be
non-negative, but this is not essential.)

To prove the proposition in the Euclidean case, consider the following
network. (See figure \ref{fig:network}).
\begin{figure}
\begin{center}
\input{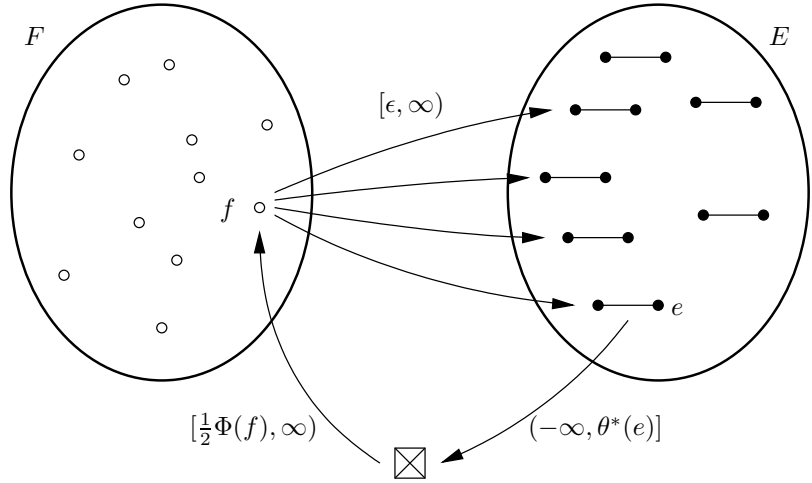}
\caption{The network $(N,X)$. 
Only a few of the branches and capacity intervals are shown.}
\label{fig:network}
\end{center}
\end{figure}
The nodes are all faces and non-oriented edges of the cellular surface, and
one further node that we denote by $\boxtimes$: $N=F\cup E\cup
\{\boxtimes\}$. There is a branch in $X$ going from $\boxtimes$ to each face
$f\in F$ with capacity interval $[\frac{1}{2}\Phi_f,\infty)$. From each face
$f$ there is a branch in $X$ going to the non-oriented edges of the boundary
of $f$ with capacity interval $[\epsilon,\infty)$, where $\epsilon>0$ will be
determined later. Finally there is a branch in $X$ going from each
non-oriented edge $e \in E$ to $\boxtimes$ with capacity
$(-\infty,\theta^*_e]$.

A feasible flow in the network yields a coherent angle system. Indeed,
Kirchoff's current law at $\boxtimes$ implies that the flow into each face
$f$ is $\frac{1}{2}\Phi_f$ and the flow out of each edge $e$ is
$\theta^*_e$. It follows that the flow in the branches from $F$ to $E$
constitutes a coherent angle system.

Assume the conditions of theorem \ref{thm:fundamental} are fulfilled.  We
need to show that the condition of the feasible flow theorem is
satisfied. Suppose $N'$ is a nonempty proper subset of $N$. Let $F'=N'\cap
F$ and $E'=N'\cap E$.

Consider first the case that $\boxtimes\in N'$, which is the easy one. Since
$N'$ is a proper subset of $N$ there is a face $f\in F$ or an edge $e\in E$
which is not in $N'$. In the first case there is a branch out of $N'$ with
infinite upper capacity bound. In the second case there is a branch into $N'$
with negative infinite lower capacity bound. Either way, the condition of the
feasible flow theorem is trivially fulfilled.

Now consider the case that $\boxtimes\not\in N'$. We may assume that for each
face $f\in F'$, the edges of the boundary of $E$ are in $E'$. Otherwise there
are branches out of $N'$ with infinite upper capacity bound. For subsets
$A,B\subset N$ denote by $A\rightarrow B$ the set of branches in $X$ having
initial node in $A$ and terminal node in $B$. Then the condition of the
feasible flow theorem is
\begin{equation*}
\sum_{f\in F'}\frac{1}{2}\Phi_f + \epsilon\,|F\setminus F'\rightarrow E'|
\leq 
\sum_{e\in E'} \theta^*_e.
\end{equation*}
It is fulfilled if we choose
\begin{equation}\label{eqn:epsilon}
\epsilon<\,\frac{1}{2|E|}
\min_{F'} \left(\sum_{e\in E'(F')}\theta^*_e
-\sum_{f\in F'}\frac{1}{2}\Phi_f\right),
\end{equation}
where the minimum is taken over all proper nonempty subsets $F'$ of $F$ and
$E'(F')$ is the set of all non-oriented edges incident with a face in
$F'$. The minimum is greater than zero because of condition (ii) of theorem
\ref{thm:fundamental}.

In the hyperbolic case, the proof is only a little bit more complicated.  In
the network, the flow in the branches going from $\boxtimes$ to a face $f$
must constrained to be exactly $\frac{1}{2}\Phi(f)$; and the capacity
interval of branches going from an edge $e$ to $\boxtimes$ has to be changed
to $(-\infty,\theta^*(e)-\epsilon]$.
\end{proof}

%---------------------------------------------------------------------

\section{Circle patterns in the sphere, Rivin's theorem, and its higher genus
analogue}\label{sec:sphereRivinHigherGenus}

In section \ref{sec:cpSphere} we adapt theorem \ref{thm:fundamental} to
circle patterns in the sphere. Then we use it to prove Rivin's theorem in
section \ref{sec:RivinProof}. Section \ref{sec:hiGenus} is devoted to the
proof of theorem \ref{thm:higherGenusRivin}, the higher genus analogue.

\subsection{Circle patterns in the sphere}\label{sec:cpSphere}

\begin{figure}
\hfill
\subfigure[{}]{\label{fig:cube}\includegraphics{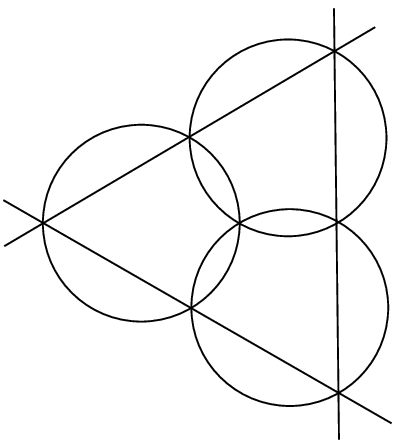}}
\hfill 
\subfigure[{}]{\label{fig:boundCirc}%
\raisebox{3mm}{\input{vacp-boundcirc.pstex_t}}}
\hspace*{\fill}
\caption{(a) The regular cubical pattern after stereographic projection. (b)
How to calculate the $\Phi_f$ for the new boundary circles.}
%\label{fig:spherePattern}
\end{figure}

A circle pattern in the sphere may be projected stereographically to the
plane, choosing some vertex, $v_{\infty}$, as the center of projection. 
For example, figure \ref{fig:cube}
shows the circle pattern combinatorially equivalent to the cube and with
intersection angles $\pi/3$ after stereographic projection. 

One obtains a circle pattern in the plane in which some circles (those
corresponding to faces incident with $v_{\infty}$) have degenerated to
straight lines. Since stereographic projection is conformal, the
intersection angles are the same. Furthermore, a M{\"o}bius-equivalent circle
pattern in the sphere leads to a planar pattern which is similar, provided
the same vertex is chosen as the center of projection.

The idea is to construct a circle pattern in the sphere by constructing the
corresponding planar pattern using the Euclidean functional and then
projecting it to the sphere. 

To construct the planar pattern, we first remove from the cellular surface
the faces incident with $v_{\infty}$ and all edges incident with them. Some
of the remaining faces will then be boundary faces. For such a new boundary
face $f$ we set 
\begin{equation}\label{eqn:PhiBoundary}
\Phi_f=2\pi-\sum 2\theta^*_e, 
\end{equation}
where the sum is taken over all removed edges incident with $f$; see figure
\ref{fig:boundCirc}. For all interior faces $f$, we set 
\begin{equation}\label{eqn:PhiInterior} 
\Phi_f=2\pi.
\end{equation}

Assuming the conditions of theorem \ref{thm:fundamental} are satisfied, we
then construct the corresponding planar pattern, add the lines
corresponding to the removed faces, and project to the sphere.

In this way, we deduce the following proposition from theorem
\ref{thm:fundamental}.

\begin{proposition}\label{prop:sphere}
Let $\Sigma$ be a closed cellular surface of genus $0$, and let an angle be
assigned to each non-oriented edge by $\theta\in(0,\pi)^E$. Suppose that for
each vertex $v\in V$, the sum of $\theta$ over the edges around $v$
is $2\pi$, so that there is no curvature in the vertices.  Let
$\theta^*=\pi-\theta$.  

Suppose $v_{\infty}\in V$, and let $F_{\infty}$ be the set of all faces
incident with $v_{\infty}$. Suppose the following conditions holds.
\begin{conditionlist}
\item{(i) } The $1$-skeleton of the Poincar{\'e} dual $\Sigma^*$
remains connected if the faces in $F_{\infty}$ and all edges incident with
them are removed.
\item{(ii) }
\begin{equation}\label{eqn:sphereEq}
2\pi(|F|-|F_{\infty}|)=\sum 2\theta^*(e),
\end{equation}
where the sum is taken over all edges incident with a face in $F\setminus
F_{\infty}$, i.~e.\, all edges not incident with $v_{\infty}$.
\item{(iii) } If $F'$ is a nonempty subset of $F\setminus
  F_{\infty}$ and $F'\not= F\setminus F_{\infty}$, and $E'$ is the set of all
  edges incident with a face in $F'$, then
\begin{equation}\label{eqn:sphereIneq}
2\pi|F'|<\sum_{e\in E'}2\theta^*(e).
\end{equation}
\end{conditionlist}
Then there exists a circle pattern in the sphere which is combinatorially
equivalent to $\Sigma$ and which has intersection angles $\theta^*$. This
circle pattern is unique up to M{\"o}bius transformations. 
\end{proposition}

\begin{proof}
By condition (i), after we remove from $\Sigma$ all faces in $F_{\infty}$ and
all edges incident with them, we are left either with a single face, or with
a connected cellular surface $\Sigma_0$ with boundary faces and boundary
vertices. The first case is elementary. In the second case, we want to apply
theorem \ref{thm:fundamental}. 

For boundary and interior faces $f$ of $\Sigma_0$ define $\Phi_f$ by
equations (\ref{eqn:PhiBoundary}) and (\ref{eqn:PhiInterior}),
respectively. We have to show that for boundary faces $\Phi_f>0$. Suppose $f$
is a boundary face of $\Sigma_0$. Let
$F'=F\setminus(F_{\infty}\cup\{f\})$. Subtracting (\ref{eqn:sphereEq}) and
(\ref{eqn:sphereIneq}), one obtains $2\pi>\sum 2\theta^*$, where the sum is
taken over all removed edges incident with $f$.  Hence $\Phi_f>0$.

We apply theorem \ref{thm:fundamental} to the cellular surface $\Sigma_0$
with intersection angles given by $\theta$ restricted to the edges of
$\Sigma_0$ and with $\Phi$ as defined above. Conditions (ii) and (iii) of the
proposition imply conditions (i) and (ii) of theorem \ref{thm:fundamental}.
Indeed, the face set of $\Sigma_0$ is $F\setminus F_{\infty}$ and
\begin{equation*}
\sum_{f\in F\setminus F_{\infty}}\Phi_f = 2\pi(|F|-|F_{\infty}|)
-\sum 2\theta^*_e,
\end{equation*}
where the sum on the right hand side is taken over all edges of $\Sigma$ that
are not edges of $\Sigma_0$ but are incident with a face in $F\setminus
F_{\infty}$. Hence, condition (ii) of the proposition implies 
\begin{equation*}
\sum_{f\in F\setminus F_{\infty}}\Phi_f = \sum 2\theta^*_e,
\end{equation*}
where the sum is taken over all edges of $\Sigma_0$, i.~e.\ condition (i)
of theorem \ref{thm:fundamental}.

Similarly, condition (iii) of the proposition implies condition (ii) of
theorem \ref{thm:fundamental}. Therefore, the corresponding planar circle
pattern exists, and by the construction described above, we can construct the
spherical circle pattern.
\end{proof}

\subsection{A proof of Rivin's theorem}\label{sec:RivinProof}

We now give a proof of Rivin's theorem.  First, we show that the condition of
Rivin's theorem is necessary. If $\gamma$ is the coboundary of a single
vertex, then $\sum_{e\in\gamma}\theta_e=2\pi$ simply means that the
corresponding circle pattern has no curvature at that vertex. Suppose a
circle pattern on the sphere is given and $\gamma$ is a closed copath
that does not bound a single vertex.  Without loss of generality, we may
assume that $\gamma$ is simple, i.~e.\ that it meets a face at most once.

\begin{figure}
\begin{center}
\input{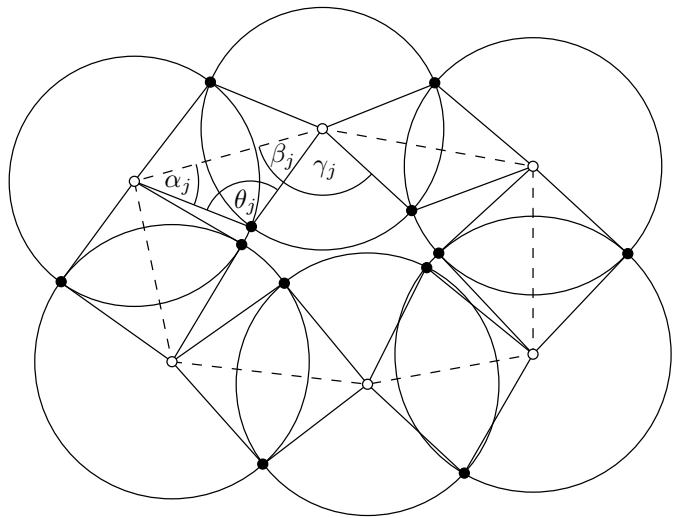}
\end{center}
\caption{The circles and quadrilaterals corresponding to a simple closed
  copath after stereographic projection.}
\label{fig:copath}
\end{figure}

Suppose that $\gamma$ does not meet all the faces of the cellular
decomposition; we will deal with the other case afterwards. Then there is a
point in the sphere which lies outside of all the circles belonging to
$\gamma$. (In fact, each circle of a circle pattern contains points which are
not contained in any of the other circles.) Project the circle pattern
stereographically to the plane, using such a point as the center of
projection. One obtains a sequence of intersecting circles as in figure
\ref{fig:copath}.
Now 
\begin{equation}\label{eqn:alphaBetaGamma}
  {\textstyle \sum} (\pi-\alpha_j-\beta_j-\gamma_j)=2\pi
\end{equation}
and
\begin{equation}\label{eqn:alphaBetaTheta}
  \theta_j=\pi-\alpha_j-\beta_j
\end{equation}
This implies $\sum\theta_j\geq 2\pi$,
with equality if and only if all $\gamma_j$ are zero. This is the case only
if the copath bounds a single vertex. 

Now suppose that the copath $\gamma$ does meet all the faces of the cellular
decomposition. Consider the graph $(V,E\setminus\gamma)$, consisting of all
vertices of the cellular decomposition and all edges that do not appear in
$\gamma$. This graph is the disjoint union of two trees. Therefore,
$|V|-|E\setminus\gamma|=2$. Since $\theta$ sums to $2\pi$ around each vertex,
\begin{equation}\label{eqn:thetaSum}
  2\pi|V|=\sum_{e\in E}2\theta_e.
\end{equation}
With $\theta_e<\pi$ one obtains,
\begin{equation*}
  \sum_{\gamma}2\theta_e=2\pi|V|-\sum_{E\setminus\gamma}2\theta_e
  >2\pi(|V|-|E\setminus\gamma|),
\end{equation*}
and hence the desired inequality. Hence the condition of Rivin's theorem is
necessary for the existence of a circle pattern.

Now, using proposition \ref{prop:sphere}, we show that the condition of
Rivin's theorem is sufficient for the existence of a circle pattern. Thus we
assume that $\sum_{e\in\gamma}\theta_e\geq 2\pi$ for every copath $\gamma$,
with equality if and only if $\gamma$ is the the coboundary of a single
vertex. We need to deduce the conditions of proposition \ref{prop:sphere}.
Condition (i) follows from the assumption that $\Sigma$ is polytopal. 

To deduce conditions (ii) and (iii), suppose that $F'$ is a subset of
$F\setminus F_{\infty}$, and $E'$ is the set of all edges of $\Sigma$
incident with a face in $F'$.  Consider $F''=F\setminus F'$, the complement
of $F'$, and $E''=E\setminus E'$, the complement of $E'$. For each edge in
$E''$, the faces on either side are in $F''$.  Now consider the dual cellular
decomposition $\Sigma^*$, and in it, the graph $\Gamma=(F'', E'')$ with
vertex set $F''$ and edge set $E''$. As for any graph, we have
\begin{equation*}
|E''|-|F''|=c-n,
\end{equation*}
where $n$ is the number of connected components of $\Gamma$ and $c$ is the
dimension of the cycle space. Since the graph is embedded in $\Sigma^*$, a
cellular decomposition of the sphere, we have 
\begin{equation*}
c=r-1,
\end{equation*}
where $r$ is the number of regions into which $\Gamma$ separates
$\Sigma^*$. Since $E''$ contains the edges of $\Sigma$ incident with
$v_{\infty}$, or, dually, the edges of $\Sigma^*$ in the boundary of
$v_{\infty}$, {\em the number of regions is at least two}. 

The boundary of each region is a copath in $\Sigma$, and hence the sum of
$\theta$ over each boundary is at least $2\pi$. Sum over all boundaries to
obtain
\begin{equation*}
2\pi r\leq\sum_{e\in E''}2\theta_e.
\end{equation*}
Indeed, each edge in $E''$ appears in at most $2$ boundaries. Equality holds
if and only if every edge is contained in a boundary {\em and} each region contains only one vertex or all vertices but one. This means equality holds
if and only if $E''=E$ or $E''=\partial^*{v}$ for some vertex
$v$; equivalently, if and only if $F'=\emptyset$ or $F'=F\setminus
F_{\infty}$.

Thus, we have shown that
\begin{equation}\label{eqn:R1}
2\pi(|E''|-|F''|)\leq \sum_{e\in E''}2\theta_e-2\pi(n+1),
\end{equation} 
with equality if and only if $F'=\emptyset$ or $F'=F\setminus F_{\infty}$. 

Equation (\ref{eqn:thetaSum}) and Euler's formula,
\begin{equation*}
|F|-|E|+|V|=2,
\end{equation*}
imply
\begin{equation}\label{eqn:R2}
2\pi(|F|-|E|)=-\sum_{e\in E}2\theta_e+4\pi.
\end{equation}
Add equation (\ref{eqn:R2}) and inequality (\ref{eqn:R1}) to obtain
\begin{equation*}
2\pi(|F'|-|E'|)\leq -\sum_{e\in E'}2\theta_e-2\pi(n-1),
\end{equation*}
or, because $\theta^*=\pi-\theta$, 
\begin{equation*}
2\pi|F'|\leq \sum_{e\in E'}2\theta^*_e+2\pi(1-n).
\end{equation*}
Equality holds if and only if $F'=\emptyset$ or $F'=F\setminus F_{\infty}$.
Since $n\geq 1$, and $n=1$ if $F'=F\setminus F_{\infty}$, we have deduced
conditions (ii) and (iii) of proposition \ref{prop:sphere}. This completes
the proof of Rivin's theorem.

\subsection{Higher genus. Proof of theorem \ref{thm:higherGenusRivin}}
\label{sec:hiGenus}

In this section, proceeding in a similar way as in the last section, we
derive theorem \ref{thm:higherGenusRivin} from theorem \ref{thm:fundamental}.

In one way, the higher genus case is simpler than the spherical case. We do
not have to transfer the problem to the plane by a stereographic projection.
We may use theorem \ref{thm:fundamental} directly. For this reason, we do not
have to assume the cellular surface to be polytopal.

In another way this case is more complicated. When a graph is embedded in
a surface of genus $g\geq 1$, the dimension $c$ of its cycle space and the
number $r$ of regions into which the graph separates the surface are not
related by the simple formula $c=r-1$. Instead, the following proposition
holds.

\begin{proposition}\label{prop:generalizedEuler}
  Suppose a non-empty graph\/ $\Gamma=(V,E)$ is embedded in a closed surface
  of genus $g$ and separates the surface into $r$ regions. Then
  \begin{equation}\label{eqn:generalizedEuler}
    r-|E|+|V| = 2-2g + \sum_{j=1}^{r} h_j,
  \end{equation}
  where $h_j$ is the dimension of the first $\mathbb{Z}_2$-homology group of
  the $j^{\text{th}}$ region.
\end{proposition}
(The graph\/ $\Gamma$ is non-empty if $V\not=\emptyset$.) If all regions are
simply connected, then all $h_j=0$, and equation (\ref{eqn:generalizedEuler})
is Euler's equation for a surface of genus $g$.  A multiply connected region
can be turned into a simply connected one by adding $h_j$ edges which do not
separate the region.  A formal proof of proposition
\ref{prop:generalizedEuler} is given by Giblin \cite[ch.~9]{Giblin1977}. It
is an application of the long exact sequence of relative homology and
Lefschetz duality.

First, we show that the condition of theorem \ref{thm:higherGenusRivin}
is necessary for the existence of a circle pattern. For $g=1$, the flat case,
this follows from equations (\ref{eqn:alphaBetaGamma}) and
(\ref{eqn:alphaBetaTheta}) as in the preceding section. (See figure
\ref{fig:copath}.) For $g>1$, one has instead
\begin{equation*}
  \textstyle\sum(\pi-\alpha_j-\beta_j-\gamma_j)=2\pi+A
\end{equation*} 
and
\begin{equation*}
  \theta_j=\pi-A_j-\alpha_j-\beta_j,
\end{equation*} 
where $A$ is the area enclosed by $\gamma$ and $A_j$ is the area of the
triangle with angles $\alpha_j$, $\beta_j$ and $\theta_j$. Hence,
$\sum\theta_j=2\pi+\sum\gamma_j+A-\sum A_j$, and therefore $\sum\theta_j\geq
2\pi$ with equality if and only if $\gamma$ is the coboundary of a single
vertex.

In the following, we show that the condition of theorem
\ref{thm:higherGenusRivin} is also sufficient for the existence of a circle
pattern.  Suppose the conditions of the theorem \ref{thm:higherGenusRivin}
are satisfied. We deduce the conditions of theorem \ref{thm:fundamental}.
Thus, let $F'$ be a non-empty subset of $F$ and let $E'$ be the set of edges
incident with any face in $F'$. Let $F''=F\setminus F'$, $E''=E\setminus E'$,
and consider the graph $\Gamma=(F'',E'')$ in the dual cellular surface
$\Sigma^*$. This graph is non-empty, if $F'\not=F$. In that case, proposition
\ref{prop:generalizedEuler} implies
\begin{equation}\label{eqn:generalizedEulerApplied}
|E''|-|F''|=2g-2+r-\sum_{j=1}^{r}h_j.
\end{equation}

Suppose we cut the dual surface $\Sigma^*$ along the (dual) edges in $E''$.
This would separate the surface into pieces. By the boundary of one such
piece, we mean the boundary after having cut the surface, not the boundary of
the corresponding region in the surface. For example, if and edge $e\in E''$
has the same region on either side, the boundary of that piece contains $e$
twice. Thus, every edge of $E''$ appears exactly twice in the boundaries of
the pieces into which the surface separates. Hence, if $\gamma_j$ is the
boundary of the $j^{\text{th}}$ region,
\begin{equation*}
  \sum_{j=1}^r\bigg(\sum_{e\in\gamma_j}\theta_e\bigg)
  =\sum_{e\in E''}2\theta_e.
\end{equation*}
By the condition of theorem \ref{thm:higherGenusRivin}, for a simply
connected piece, i.~e.\ for one with $h_j=0$, the sum of $\theta$ over its
boundary is at least $2\pi$, and equal to $2\pi$ only if the piece consists
of a single vertex. Hence, we have
\begin{equation*}
  2\pi(r-\sum_{j=1}^{r}h_j)<\sum_{e\in E''}2\theta_e.
\end{equation*}
Equality could hold only if each region consisted of a single vertex. But
this would imply $E''=E$, or $E'=\emptyset$, contradicting the assumption
that $F'$ was non-empty. If $F''$ is non-empty, equation
(\ref{eqn:generalizedEulerApplied}) implies
\begin{equation}\label{eqn:HG1}
2\pi(|E''|-|F''|)<2\pi(2g-2)+\sum_{e\in E''}2\theta_e.
\end{equation}
If $F''$ is empty (and hence also $E''$), this inequality is trivially true
if $g>1$, and equality holds if $g=1$.

Because there is no curvature in the vertices, equation (\ref{eqn:thetaSum})
holds here as well. With Euler's formula,
\begin{equation*}
 |F|-|E|+|V|=2-2g,
\end{equation*}
this implies
\begin{equation*}%\label{eqn:HG2}
2\pi(|F|-|E|)=2\pi(2-2g)-\sum_{e\in E}2\theta_e.
\end{equation*}
Add this equation to inequality (\ref{eqn:HG1}) to obtain
\begin{equation*}
2\pi(|F'|-|E'|)\leq-\sum_{e\in E'}2\theta_e,
\end{equation*}
with equality if and only if $F'=F$ and $g=1$. Since $\theta^*=\pi-\theta$,
this inequality is equivalent to 
\begin{equation*}
2\pi|F'|\leq\sum_{e\in E'}2\theta_e^*.
\end{equation*}
This completes the proof of theorem \ref{thm:higherGenusRivin}.

%---------------------------------------------------------------------

\section{Other variational principles}\label{sec:other}

In this section, we derive the variational principles of
Colin~de~Verdi{\`e}re, Br{\"a}gger, and Rivin from ours.  For circle
packings, Colin~de~Verdi{\`e}re's functionals are obtained from our
functionals $\Seuc$ and $\Shyp$ for orthogonally intersecting circles by
minimizing with respect to the radii of those circles which do not appear in
the pattern.  The derivation of Br{\"a}gger's and Rivin's functionals
involves a Legendre transformation of $\Seuc$. We believe that Leibon's
functional can be derived in the same way from $\Shyp$ by the Legendre
transformation which we present.

\subsection{Colin~de~Verdi{\`e}re's functionals}\label{sec:CdV}

\begin{figure}
\hfill
\input{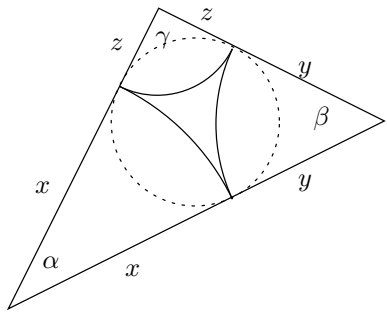}
\hspace*{\fill}
\caption{A Euclidean triangle}
\label{fig:eucTri}
\end{figure}

Colin~de~Verdi{\`e}re considers circle packings in which the circles
correspond to the vertices of a triangulation. He considers the $1$-form 
\begin{equation*}
\omega = \alpha\,du+\beta\,dv+\gamma\,dw
\end{equation*}
on the space of Euclidean triangles, where 
\begin{equation*}
  u=\log x, \quad v=\log y, \quad\text{and}\quad w=\log z
\end{equation*}
and $x$, $y$, $z$ and $\alpha$, $\beta$, $\gamma$ are as shown in figure
\ref{fig:eucTri}.

It turns out that $d\omega=0$, hence one may integrate. Define the
function $f_{\alpha_0, \beta_0, \gamma_0}$ on $\mathbb{R}^3$ by
\begin{equation*}
  f_{\alpha_0, \beta_0, \gamma_0}(u,v,w) = 
  \int^{(u,v,w)} (\alpha_0 - \alpha)\,du + (\beta_0 - \beta)\,dv
  + (\gamma_0 - \gamma)\,dw.
\end{equation*}
The initial point of the integration does not matter. 

Suppose we are given a triangulation and a coherent angle system for it.
Here, a coherent angle system is a positive function on the set of angles of
the triangles, such that the sum in each triangle is $\pi$, and the sum
around each vertex is $2\pi$. For a function $\rho$ on the vertices of the
triangulation, Colin~de~Verdi{\`e}re's functional for Euclidean circle
packings is
\begin{equation*}
  S_{\text{\it CdV}}(\rho)=\sum f_{\alpha_0, \beta_0, \gamma_0}(u, v, w),
\end{equation*}
where the sum is taken over all triangles, $u$, $v$, and $w$ are the values
of $\rho$ on the vertices of each triangle, and $\alpha_0$, $\beta_0$ and
$\gamma_0$ are the corresponding angles of the coherent angle system.

The critical points of this functional correspond to the logarithmic radii of
a circle packing.

The hyperbolic case is treated in the same way, except that now 
\begin{equation*}
    u=\log \tanh(x/2), \quad v=\log \tanh(y/2), 
    \quad\text{and}\quad w=\log \tanh(z/2),
\end{equation*}
and in the definition of a coherent angle system, it is required that the
sum of the angles in a triangle is less than $\pi$.

To treat circle packings with our functionals, we consider circle patterns
with orthogonally intersecting circles on the medial decomposition of the
triangulation (see figure \ref{fig:medial}). 

Suppose $\Sigma_0$ is a triangulation. Let $\Sigma$ be the medial
decomposition. Set $\theta_e=\pi/2$ for all edges $e$. We will consider only
the case of closed surfaces, so let $\Phi_f=2\pi$ for all faces $f$. The
faces of $\Sigma$ are of two types: those that correspond to faces of
$\Sigma_0$ and those that correspond to vertices of $\Sigma_0$.  Let
$F=F_1\cup F_2$, where $F_1$ contains the faces of the first type and $F_2$
contains the faces of the second type.

Consider $\Seuc$ as function on $\mathbb{R}^{F_1}\times\mathbb{R}^{F_2}$ and
define $S_2:\mathbb{R}^{F_2}\rightarrow \mathbb{R}$,
\begin{equation*}
S_2(\rho_2)=\min_{\rho_1\in\mathbb{R}^{F_1}} 
\Seuc(\rho_1,\rho_2).
\end{equation*}

Hence, $S_2(\rho_2)=\Seuc(\rho_1,\rho_2)$, where $\rho_1\in\mathbb{R}^{F_1}$
is determined as follows.  Suppose $f\in F_1$ and $f_a, f_b, f_c\in F_2$ are
the neighboring faces of $f$. Construct the Euclidean triangle whose sides
are
\begin{equation*}
\e{{\rho_2}_{f_a}}+\e{{\rho_2}_{f_b}},\quad
\e{{\rho_2}_{f_b}}+\e{{\rho_2}_{f_c}},\quad\text{and}\quad
\e{{\rho_2}_{f_c}}+\e{{\rho_2}_{f_a}}.
\end{equation*}
Let $\rho_1(f)$ be the logarithmic radius of the inscribed circle.

For $\rho_2\in\mathbb{R}^{F_2}$, and $\rho_1$ the corresponding
point in $\mathbb{R}^{F_1}$,
\begin{equation*}
  \frac{\partial \Seuc}{\partial \rho_1}(\rho_1, \rho_2)  = 0,\quad
  \frac{\partial \Seuc}{\partial \rho_2}(\rho_1, \rho_2) =
  \frac{\partial S_2}{\partial \rho_2}(\rho_2).
\end{equation*}
Consider $S_{\text{\it CdV}}$ as function on $\mathbb{R}^{F_2}$. It is not
hard to see that $dS_2=\,dS_{\text{\it CdV}}$.  This implies that
Colin~de~Verdi{\`e}re's Euclidean functional is, up to an additive constant,
equal to $S_2$.

Colin~de~Verdi{\`e}re's functional for hyperbolic circle packings
can be derived from $\Shyp$ in the same way.

\subsection{A Legendre transformation of the Euclidean functional}

Consider a function $S:\mathbb{R}^F\rightarrow\mathbb{R}$ of the form
\begin{equation*}
  S(\rho) = \sum_{f_j\circ\edgeup\circ f_k} F_e(\rho_{f_k}-\rho_{f_j})
  + \sum_{\circ f_j} G_f(\rho_{f_j}),
\end{equation*}
where the first sum is taken over all oriented edges $e$ and the second sum
is taken over all faces $f$. The functional $\Seuc$ is of this form with 
\begin{equation*}
  F_e(x)=\im\Li\big(\e{x+i\theta_e})
\end{equation*}
and
\begin{equation*}
  G_f(x)=\Big(\Phi_f - \sum \theta^*_e\Big)x,
\end{equation*}
where the sum is taken over all edges in the boundary of $f$.

Define the function
$S^L:\mathbb{R}^{\vec{E}}\times\mathbb{R}^{F}\rightarrow\mathbb{R}$ by
\begin{equation*}
  S^L(v, \rho) = \sum_{f_j\circ\edgeup\circ f_k} F_e(v_e)
  + \sum_{\circ f_j} G_{f_j}(\rho_{f_j}).
\end{equation*}
The `L' in the superscript stands for `Lagrange'. We do not assume
that $v$ is antisymmetric, i.~e.\ that $v_e=v_{-e}$. Clearly, the critical
points of $S$ correspond to the critical points of $S^L$ under the
constraints
\begin{equation*}
  v_e=\rho_{f_k}-\rho_{f_j},
\end{equation*}
where $f_j$ and $f_k$ are the faces to the left and to the right of $e$,
respectively.

\begin{definition}
  Suppose $F$ is a smooth and strictly convex real function on some open
  interval. Then
  \begin{equation}\label{eqn:yOfx}
    y=F'(x)
  \end{equation}
  defines a smooth coordinate transformation. The {\em Legendre transform} of
  $F$ is
  \begin{equation*}
    \tilde{F}(y) = xy - F(x),
  \end{equation*}
  where $x$ is related to $y$ by (\ref{eqn:yOfx}).
\end{definition}

Note that
\begin{equation*}
\tilde{F}'(y)=\frac{d\tilde{F}}{dy}= \frac{dx}{dy} y + x 
- \frac{dF}{dx}\frac{dx}{dy} = x.
\end{equation*}
It follows that $F$ is the Legendre transform of $\tilde{F}$. 

Now let $\tilde{F}_e$ be the Legendre transform of $F_e$ and define the
function $S^H$ by
\begin{equation*}
    S^H(\varphi, \rho) = \sum_{f_j\circ\edgeup\circ f_k} 
    \big(\varphi_e(\rho_{f_k}-\rho_{f_j})-\tilde{F}_e(\varphi_e)\big)
  + \sum_{\circ j} G(\rho_{f_j}).
\end{equation*}
The $H$ in the superscript stands for `Hamilton'. If $\varphi$ and $v$ are
related by $\varphi_e=F_e'(v_e)$, then $S^H(\varphi, \rho)=S^L(v,
\rho)$. Remarkably, the critical points of $S$ correspond to critical points
of $S^H$, where $\rho$ and $\varphi$ may be varied independently:
\begin{proposition}\label{prop:criticalPoints}
  If $(\varphi, \rho)$ is a critical point of $S^H$, then $\rho$ is a
  critical point of $S$. 
  
  Conversely, suppose $\rho$ is a critical point of $S$. Let
  $\varphi_e=F_e'(\rho_{f_k}-\rho_{f_j})$, where $f_j$ and $f_k$ are the
  faces to the left and right of $e$. Then $(\varphi, \rho)$ is a critical
  point of $S^H$. 
\end{proposition}
\begin{proof}
  We have $\partial S^H/\partial\varphi_e =
  \rho_{f_k}-\rho_{f_k}-\tilde{F}_e'(\varphi_e)$.  So all partial derivatives
  with respect to the variables $\varphi_e$ vanish, if and only if $\varphi$
  and $\rho$ are related by $\rho_{f_k}-\rho_{f_k}=\tilde{F}_e'(\varphi_e)$,
  or, equivalently, $\varphi_e=F_e'(\rho_{f_k}-\rho_{f_k})$. In this case
  $S^H(\varphi, \rho)=S(\rho)$.
\end{proof}

Now consider in particular the Euclidean functional $\Seuc$. The variables
$\rho$ and $\varphi$ are related by equation (\ref{eqn:phiOfRhoEuc}). The
inverse relation is 
\begin{equation}\label{eqn:rhoOfPhiEuc}
  \rho_{f_k}-\rho_{f_j}=\log\frac{\sin\varphi_e}{\sin(\varphi_e+\theta_e)}.
\end{equation}
For the functional in Hamiltonian form, identity (\ref{eqn:LiCl}) of the
appendix implies
\begin{multline*}
  \SeucH(\varphi, \rho)=
  \sum_{f_j\circ\edgeup\circ f_k}
  \Big(\varphi_e\,(\rho_{f_k}-\rho_{f_j})
  +\frac{1}{2}\big(\Cl(2\varphi_e)-\Cl(2\varphi_e+2\theta_e)
  +\Cl(2\theta_e) \big)\Big)\\
  + \sum_{\circ f}
  \Big(\Phi_f - \sum_{e\in\partial f} \theta^*_e \Big) \rho_f.
\end{multline*}

Note that $\SeucH$ depends linearly on $\rho$. Collecting the coefficients of
each $\rho_f$, one finds that if $\varphi$ is a coherent angle system, then
$\SeucH(\varphi, \rho)$ does not depend on $\rho$ at all. Therefore, this
leads to a new variational principle, where the functional depends only on
the angles, and the variation is constrained to the space of coherent angle
systems.
\begin{proposition}\label{prop:SeucHred}
  If $\varphi$ is a coherent angle system, then $\SeucH(\varphi,
  \rho)=\SeucHred(\varphi)$, where
  \begin{equation}\label{eqn:SeucHred}
    \SeucHred(\varphi)
    =\sum \Big(\Cl(2\varphi_e)+\frac{1}{2}\Cl(2\theta_e)\Big)
  \end{equation}
  and the sum is taken over all oriented edges $e$.
  
  If $\varphi$ is a critical point of $\SeucHred$ under variations in the
  space of coherent angle systems, then the equations (\ref{eqn:rhoOfPhiEuc})
  form a consistent system of linear equations for the $\rho_f$. They define
  a function $\rho\in\mathbb{R}^F$ up to an additive constant. This $\rho$ is
  a critical point of $\Seuc$.
  
  Conversely, any critical point of $\Seuc$ yields, via equation
  (\ref{eqn:phiOfRhoEuc}), a coherent angle system which is a critical point
  of $\SeucHred$ under variations in the space of coherent angle systems.
\end{proposition}
\begin{proof}
  To derive equation (\ref{eqn:SeucHred}), note that, if $\varphi$ is a
  coherent angle system, then $\varphi_{e}+\theta_e=\pi-\varphi_{-e}$. Hence
  $\Cl(2\varphi_e+2\theta_e)=-\Cl(2\varphi_{-e})$.
  
  Now suppose $\varphi$ is a critical point of $\SeucHred$ under variations
  in the space of coherent angle systems. The partial derivatives of
  $\SeucHred$ are 
  \begin{equation*}
    \frac{\partial\SeucHred}{\partial\varphi_e}=-2\log(2\sin\varphi_e).
  \end{equation*}
  The tangent space to the space of coherent angle systems is spanned by
  vectors of the form
  \begin{equation*}
    \sum_{e\in\gamma}\Big(\frac{\partial}{\partial\varphi_e}
    -\frac{\partial}{\partial\varphi_{-e}}\Big),
  \end{equation*}
  where $\gamma$ is some cocycle. (If we add to $\varphi_e$ we have to
  subtract the same amount from $\varphi_{-e}$. But the sum of $\varphi$
  around each face has to remain constant, so we step from one face to the
  next, adding and subtracting along some cocycle.) 
  
  Since $\varphi$ is a critical point under variations in the space of
  coherent angle systems,
  \begin{equation*}
    \sum_{e\in\gamma}\log\frac{\sin\varphi_e}{\sin(\varphi_e+\theta_e)}
    =0
  \end{equation*}
  for all cocycles $\gamma$. Hence, equations (\ref{eqn:rhoOfPhiEuc}) are
  compatible.  It is now easy to complete the proof.
\end{proof}

\subsection{Br{\"a}gger's functional and Rivin's functional}

Like Colin~de~Verdi{\`e}re, Br{\"a}gger \cite{Bragger1992} considers circle
packings in which the circles correspond to the vertices of a triangulation.
A coherent angle system in the sense of section \ref{sec:CdV} is clearly
equivalent to a coherent angle system in the sense of this article. Thus,
Br{\"a}gger's functional is seen to be equal to $\SeucHred/2$ up to an
additive constant.

Rivin \cite{Rivin1994} considers Euclidean circle patterns with arbitrary
prescribed intersection angles. The pattern of intersection is determined by
a triangulation. Circles correspond to faces of the triangulation. However,
since Rivin allows intersection angles $\theta=0$, adjacent triangles may in
effect belong to the same circle. Thus, cellular decompositions with
non-triangular cases may be treated by first dissecting all faces into
triangles. Rivin treats cone-like singularities in the vertices but not in
the centers of the circles. His functional is up to an additive constant
equal to $\SeucHred/2$.

\subsection{A Legendre transformation of the hyperbolic Functional}

Now, consider a function $S:\mathbb{R}^F\rightarrow\mathbb{R}$ of the form
\begin{equation*}
  S(\rho) = \sum_{f_j\circ\edge\circ f_k} 
  \big(F^-_e(\rho_{f_k}-\rho_{f_j})
  + F^+_e(\rho_{f_j}+\rho_{f_k})\big)
  + \sum_{\circ f} G_f(\rho_f),
\end{equation*}
Here, the first sum is taken over all {\em non}\/oriented edges $e$, and
$f_j$ and $f_k$ are the faces on either side of $e$. We assume that $F^-$ is
an even function so that is does not matter which face is which. The second
sum is taken over all faces $f$. The functional $\Shyp$ is of this form with
\begin{equation*}
  F^-_e(x) = F^+_e(x) 
  = \im\Li\big(\e{x+i\theta_e}) + \im\Li\big(\e{-x+i\theta_e})
\end{equation*}
and
\begin{equation*}
  G_f(x)=\Phi_f\, x.
\end{equation*}

The space $\mathbb{R}^{\vec{E}}$ of functions on non-oriented edges splits
into the space of antisymmetric and symmetric functions,
\begin{equation*}
  \mathbb{R}^{\vec{E}}=\Alt(\vec{E})\oplus\Sym(\vec{E}),
\end{equation*}
where 
\begin{equation*}
  \Alt(\vec{E})=\big\{v\in\mathbb{R}^{\vec{E}}\,\big|\,v_{-e}=-v_{e}\big\}
\end{equation*}
and
\begin{equation*}
  \Sym(\vec{E})=\big\{w\in\mathbb{R}^{\vec{E}}\,\big|\,w_{-e}=w_{e}\big\}.
\end{equation*}
Define the function 
$S^L:\Alt(\vec{E})\times\Sym(\vec{E})\times\mathbb{R}^F\rightarrow\mathbb{R}$
by
\begin{equation*}
  S^L(v, w, \rho) = \sum_{f_j\circ\edge\circ f_k} 
  \big(F^-_e(v) + F^+_e(w)\big) + \sum_{\circ f} G_f(\rho_f).
\end{equation*}

Suppose $F^-$ and $F^+$ are smooth and strictly convex. The coordinates
conjugate to $v_e$ and $w_e$ are $p_e=F^{-\prime}(v_e)$ and
$s_e=F^{+\prime}(w_e)$. Clearly, $p_{-e}=-p_e$ and $s_{-e}=s_e$. Let
$\tilde{F}^-$ and $\tilde{F}^+$ be the Legendre transforms of $F^-$ and
$F^+$. For $p$ and $s$ in suitable subsets of $\Alt(\vec{E})$ and
$\Sym(\vec{E})$, define
\begin{equation*}
\begin{split}
  S^H(p, s, \rho) = 
  \sum_{f_j\circ\edge\circ f_k} &
  \big(p_e (\rho_{f_k}-\rho_{f_j}) - \tilde{F}^-_e(p) 
     + s_e (\rho_{f_k}+\rho_{f_j}) - \tilde{F}^+_e(s) \big) \\
  + \sum_{\circ f} & G_f(\rho_f).
\end{split}
\end{equation*}
If $p$ and $s$ are related to $\rho$ by
\begin{equation}\label{eqn:psOfRhoGeneral}
  p_e={F_e^{-\prime}}(\rho_{f_k}-\rho_{f_j}),\quad
  s_e={F_e^{+\prime}}(\rho_{f_k}+\rho_{f_j}), 
\end{equation}
where $f_j$ and $f_k$ are the face to the left and right of $e$, then
$S^H(p,s,\rho)=S(\rho)$.  Again, the critical points of $S$ correspond to
critical points of $S^H$, where $p$, $s$, and $\rho$ vary independently. The
following proposition is proved in the same way as proposition
\ref{prop:criticalPoints}.
\begin{proposition}
  If $(p, s, \rho)$ is a critical point of $S^H$, subject only to the
  constraints $p_{-e}=-p_e$ and $s_{-e}=s_e$, then $\rho$ is a
  critical point of $S$.
  
  Conversely, suppose $\rho$ is a critical point of $S$. Define $p$ and $s$
  by equations (\ref{eqn:psOfRhoGeneral}). Then $(p,s,\rho)$ is a critical
  point of $S^H$.
\end{proposition}

Now consider the functional $\Shyp$. The variables $p$ and $s$ are related to
$\rho$ by equations (\ref{eqn:pOfRho}) and (\ref{eqn:sOfRho}).
The Hamiltonian form of the functional follows from equation
(\ref{eqn:LiClp}) of the appendix:
\begin{multline}\label{eqn:ShypH}
  \ShypH(p,s,\rho)=
  \sum_{f_j\circ\edge\circ f_k} \Big( 
  p_e\big(\rho_{f_k}-\rho_{f_j}\big)
  +\Cl\big(\theta^*_e+p_e\big)+\Cl\big(\theta^*_e-p_e\big)\\
  +s_e\big(\rho_{f_j}+\rho_{f_k}\big)
  +\Cl\big(\theta^*_e+s_e\big)+\Cl\big(\theta^*_e-s_e\big)\\
  -2\Cl\big(2\theta^*_e\big) \Big)
  +\sum_{\circ f} \Phi_f\rho_f.
\end{multline}
(See equation (\ref{eqn:LiClp}) of the appendix.)

Introduce the variables 
\begin{equation*}
\varphi_e=\frac{1}{2}(p_e-s_e)
\end{equation*}
on oriented edges instead of $(p,s)$:
\begin{equation*}
\ShypH(\varphi,\rho)=\ShypH(p,s,\rho).
\end{equation*}
This functional is linear in $\rho$. As in the Euclidean case, collecting the
coefficients of each $\rho_f$, one observes that at a critical point
$\varphi$ is a coherent angle system. Also, if $\varphi$ is a coherent angle
system, then $\ShypH$ does not depend on $\rho$. (Remember that in the
hyperbolic case, coherent angle systems satisfy
$\varphi_e+\varphi_{-e}<\theta^*$ instead of
$\varphi_e+\varphi_{-e}=\theta^*$.) Again, this leads to a different
variational principle.
\begin{proposition}
  If $\varphi$ is a coherent angle system, then $\ShypH=\ShypHred$, where
  \begin{multline*}
    \ShypHred(\varphi)=
    \sum \big(
       \Cl(\theta^*_e+\varphi_e-\varphi_{-e}) 
       + \Cl(\theta^*_e-\varphi_e+\varphi_{-e})\\ 
       + \Cl(\theta^*_e+\varphi_e+\varphi_{-e}) 
       + \Cl(\theta^*_e-\varphi_e-\varphi_{-e})
              - 2\Cl(2\theta^*_e) \big).
  \end{multline*}
  The sum is taken over non-oriented edges, or more precisely, over one
  representative for each pair of oppositely oriented edges.
  
  If $\varphi$ is a critical point of $\ShypHred$ under variations in the
  space of coherent angle systems, then the equations 
  \begin{equation}\label{eqn:rhoOfPhiHyp}
    \rho_{f}=\frac{1}{2}\log
    \frac{\sin\Big(\frac{\displaystyle \theta^*-\varphi_e-\varphi_{-e}}
                   {\displaystyle 2}\Big)
          \sin\Big(\frac{\displaystyle \theta^*-\varphi_e+\varphi_{-e}}
                   {\displaystyle 2}\Big)}
         {\sin\Big(\frac{\displaystyle \theta^*+\varphi_e+\varphi_{-e}}
                   {\displaystyle 2}\Big)
          \sin\Big(\frac{\displaystyle \theta^*+\varphi_e-\varphi_{-e}}
                   {\displaystyle 2}\Big)},
  \end{equation}
  where $f$ is the face to the left of $e$, coherently define
  $\rho\in\mathbb{R}^F$. This $\rho$ is a critical point of $\Shyp$.
  
  Conversely, suppose $\rho$ is a critical point of $\Shyp$. Then equations
  (\ref{eqn:phiOfRhoHyp}) define a coherent angle system, which is a critical
  point of $\SeucHred$ under variations in the space of coherent angle
  systems.
\end{proposition}
The proof is analogous to the proof of proposition \ref{prop:SeucHred}. We
would like to remark only on two points.

First, equations (\ref{eqn:rhoOfPhiHyp}) are equivalent to equations
(\ref{eqn:phiOfRhoHyp}). This follows by a straightforward calculation.

Second, to show that equations (\ref{eqn:rhoOfPhiHyp}) are consistent if
$\varphi$ is a critical point of $\ShypHred$ under variations in the space of
coherent angle systems, proceed as follows.  Suppose $f$ is the face to the
left of $e$, and $e'=\sigma e$ is the next edge in the boundary of $f$. It
suffices to show that the value for $\rho_f$ obtained from the equation
involving $\varphi_{\pm e}$ is the same as the value obtained from the
equation involving $\varphi_{\pm e'}$. The tangential space to the space of
coherent angle systems is spanned by variations with $\dot{\varphi}_e=-1$,
$\dot{\varphi}_{e'}=1$ for two oriented edges $e$ and $e'=\sigma e$, and
$\dot\varphi=0$ otherwise. But $\partial\ShypHred/\partial\varphi_{e}$ is
twice the right hand side of equation (\ref{eqn:rhoOfPhiHyp}). This implies
the consistency of equations (\ref{eqn:rhoOfPhiHyp})

\subsection{Leibon's functional}

Like Rivin, Leibon \cite{Leibon2001} considers circle patterns with arbitrary
prescribed intersection angles. The circles correspond to the faces of a
triangulation.  Whereas Rivin treats the Euclidean case, Leibon treats the
hyperbolic case. We believe that his functional can be identified with
$\ShypHred$. 

%---------------------------------------------------------------------

\appendix

\section*{Appendix: The dilogarithm function and Clausen's integral}
%\label{app:dilog}

In this appendix, we collect everything about the dilogarithm and Clausen's
integral that is relevant for this paper. A more thorough treatment and an
extensive bibliography are contained in Lewin's monograph \cite{Lewin1981}. 

For $|z|\leq 1$, the dilogarithm function is defined by the power
series
\begin{equation*}
\Li(z)=\frac{z}{1^2}+\frac{z^2}{2^2}+\frac{z^3}{3^2}+\ldots\,.
\end{equation*}
For $|z|<1$,
\begin{equation*}
-\log(1-z)=\frac{z}{1}+\frac{z^2}{2}+\frac{z^3}{3}+\ldots,
\end{equation*}
and hence
\begin{equation}\label{dilogIntegral}
\Li(z)=-\int_{0}^{z}\frac{\log(1-\zeta)}{\zeta}\,d\zeta\,.
\end{equation}

In the light of this integral representation, one sees that the dilogarithm
can be analytically continued to the complex plane cut from $1$ to $\infty$
along the positive real axis. 

Clausen's integral $\Cl(x)$ can be defined by the imaginary part of
the dilogarithm on the unit circle:
\begin{equation*}
\begin{split}
\Cl(x)&=\im \Li(\e{ix}) \\
&=\frac{1}{2i}\big(\Li(\e{ix})-\Li(\e{-ix})\big).
\end{split}
\end{equation*}
Clausen used the integral representation we will derive below to define this
function, hence the name.  We consider Clausen's integral as a real
valued function of a real variable. It is $2\pi$-periodic and odd.  The power
series representation of the dilogarithm yields
%\begin{equation*}
%\Li(\e{ix})=\sum_{n=1}^{\infty} \frac{\cos(nx)}{n^2} +
%i \sum_{n=1}^{\infty} \frac{\sin(nx)}{n^2}
%\end{equation*}
%and hence 
the Fourier series representation for Clausen's integral,
\begin{equation*}
\Cl(x)=\sum_{n=1}^{\infty} \frac{\sin(nx)}{n^2}\,.
\end{equation*}
Substitute $\zeta=\e{i\xi}$ in the integral representation of the dilogarithm
(\ref{dilogIntegral}) to obtain, for $0\leq x \leq 2\pi$,
\begin{equation*}
\Cl(x)=-\int_{0}^{x}\log\left(2 \sin\frac{\xi}{2}\right)\,d\xi\,.
\end{equation*}
Clausen's integral is almost the same as Milnor's Lobachevski function
\begin{equation*}
\operatorname{\mbox{\cyr L}}(x)=-\int_{0}^{x}\log|2 \sin \xi|\,d\xi
%\text{\cyr L}(x)=-\int_{0}^{x}\log|2 \sin \xi|\,d\xi
=\frac{1}{2}\Cl(2x)\,.
\end{equation*}
Another related function appearing in the literature, also called
Lobachevski function, is
\begin{equation*}
\operatorname{L}(x)=-\int_{0}^{x}\log(\cos \xi)\,d\xi
=\frac{1}{2}\Cl(2x-\pi)+x \log 2\,.
\end{equation*}

We will now derive a formula expressing the imaginary part of the dilogarithm
in terms of Clausen's integral not only on the unit circle, but anywhere in
the complex plane. Suppose $x$ is real and $0<\theta<2\pi$. Substitute
$\zeta=\e{\xi+i\theta}$  in (\ref{dilogIntegral}), to obtain
\begin{equation}\label{eqn:ImLiIntegral}
\begin{split}
\im \Li(\e{x+i\theta})&=
\frac{1}{2i}\big(\Li(\e{x+i\theta})
-\Li(\e{x-i\theta})\big) \\
&=\frac{1}{2i}\int_{-\infty}^{x}
\log\left(\frac{1-\e{\xi-i\theta}}{1-\e{\xi+i\theta}}\right)\,d\xi\,.
\end{split}
\end{equation}
Now substitute 
\begin{equation*}
\eta=\frac{1}{2i}
\log\left(\frac{1-\e{\xi-i\theta}}{1-\e{\xi+i\theta}}\right)\,,
\end{equation*}
and note that inversely
\begin{equation*}
\xi=\log\left(\frac{2\sin\eta}{2\sin(\eta+\theta)}\right)\,.
\end{equation*}
Finally, one obtains
\begin{equation}\label{eqn:LiCl}
\im \Li(\e{x+i\theta})
=yx + \frac{1}{2}\Cl(2y) - \frac{1}{2}\Cl(2y+2\theta)
+ \frac{1}{2}\Cl(2\theta)\,,
\end{equation}
where
\begin{equation*}
y=\frac{1}{2i}
\log\left(\frac{1-\e{x-i\theta}}{1-\e{x+i\theta}}\right)\,.
\end{equation*}
From this, we derive the formula
\begin{multline}\label{eqn:LiClp}
\im \Li(\e{x+i\theta}) + \im \Li(\e{-x+i\theta})= \\
px + \Cl(p+\theta^*) 
+ \Cl(-p+\theta^*) - \Cl(2\theta^*)\,,
\end{multline}
where $\theta^*=\pi-\theta$, and 
\begin{equation*}
p=\frac{1}{2i}\log
\frac{(1+\e{x+i\theta^*})(1+\e{-x-i\theta^*})}
{(1+\e{x-i\theta^*})(1+\e{-x+i\theta^*})}\,.
\end{equation*}
Finally, note that $p$ and $x$ are related by
\begin{equation*}
\tan\Big(\frac{p}{2}\Big)=
\tanh\Big(\frac{x}{2}\Big)\tan\Big(\frac{\theta^*}{2}\Big)\,.
\end{equation*}

%%%%%%%%%%%%%%%%%%%%%%%%%%%%%%%%%%%%%%%%%%%%%%%%%%%%%%%%%%%%%%%%%%%%%%%%%%

%\nocite{*}

%\bibliographystyle{plain}
\bibliographystyle{alpha}

\bibliography{vacp}

%\listoffigures

\end{document}